\documentclass[preprint,10pt]{elsarticle}

\usepackage{amsmath,amssymb,amsthm}
\usepackage{enumerate}
\usepackage{cite}
\usepackage{graphicx} 
\usepackage{dsfont}

\numberwithin{equation}{section}

\newcommand{\N}{\mathbb{N}}
\newcommand{\Z}{\mathbb{Z}}
\newcommand{\R}{\mathbb{R}}
\newcommand{\PR}{\mathbb{P}}

\newcommand{\sas}{\mathcal{S}\alpha\mathcal{S}}

\newcommand{\one}{\mathds{1}}

\newtheorem{thm}{Theorem}[section]
\newtheorem{lemma}{Lemma}[section]
\newtheorem{Prop}{Proposition}[section]

\def\ce{{\cal C}}

\newcommand{\Za}[1]{\mathrm{Z}_{\alpha} \left( {#1} \right) }
\def\o{\omega}
\def\eps{\epsilon}
\def\al{\alpha}
\def\ga{\gamma}
\def\th{\theta}
\def\de{\delta}
\def\lajk{\zeta_{j,k}}
\def\lajnk{\zeta_{j,-k}}
\def\lanjnk{\zeta_{-j,-k}}
\def\ljk{\lambda_{j,k}}
\def\ljnk{\lambda_{j,-k}}
\def\lnjnk{\lambda_{-j,-k}}
\def\Th{\Theta}

\biboptions{comma,sort}

\begin{document}

\begin{frontmatter}



\title{Linear Multifractional Stable Motion: representation via Haar basis}


\author{Julien Hamonier}
\address{Laboratoire de Physique, CNRS UMR 5672, ENS Lyon\\
46, all\'ee d'Italie, F-69364 Lyon cedex, France.
\ead{julien.hamonier@ens-lyon.fr}}

\begin{abstract}
The aim of this paper is to give a wavelet series representation of Linear Multifractional Stable Motion (LMSM in brief), which is more explicit than that introduced in \citep{ayache2012lmsm}. Instead of using Daubechies wavelet, which is not given by a closed form, we use the Haar wavelet.
In order to obtain this new representation, we introduce a Haar expansion of the high and low frequency parts of the $\sas$ random field $X$ generating LMSM. Then, by using Abel transforms, we show that these series are convergent, almost surely, in the space of continuous functions.  Finally, we determine their almost sure rates of convergence in the latter space. Note that these representations of the high and low frequency parts of $X$, provide a new method for simulating the high and low frequency parts of LMSM.
\end{abstract}

\begin{keyword}
Approximation of processes \sep Linear Fractional and Multifractional Stable Motions \sep Wavelet series representations \sep  Haar system.
\MSC[2010] 60G22 \sep 60G52 \sep 41A30 \sep 41A58
\end{keyword}

\end{frontmatter}


\section{Introduction and main results}

For several years, there is a growing interest in probabilistic models based on fractional and multifractional processes; they are a convenient tools for applications in various areas such as modelling  of Internet traffic, finance, \ldots These models are natural extensions of the well-known Gaussian Fractional Brownian Motion (FBM in brief), which has stationary increments and is self-similar with self-similar exponent $H\in (0,1)$. One of the most known extensions of FBM in the setting of heavy-tailed stable distributions is Linear Fractional Stable Motion (LFSM in brief) (see \citep{SamTaq,EmMa}). It has also stationary increments and is self-similar with self-similar exponent $H\in (0,1)$; nevertheless it depend on a second parameter, denoted by $\al \in (0,2)$ which control the tail heaviness of the distribution of LFSM.  In order to overcome some limitations due to the stationarity of increments and the constancy of the exponent of self-similarity of each of these processes, Linear Multifractional Stable Motion (LMSM in brief) was introduced in \citep{stoev2005path,stoev2004stochastic}; according to these two authors, LMSM is a good candidate to describe some features of Internet traffic, for example, burstiness, that is the presence of rare but extremely busy periods of activity.

To precisely define LMSM, we need to fix some notations to be used throughout the article:  
\begin{enumerate}
\item We assume that $\al\in (1,2)$, since Stoev and Taqqu, in \citep{stoev2004stochastic}, showed that this assumption is a necessary condition for the continuity, with probability 1, of the LMSM's paths.
\item $H(\cdot)$ denotes an arbitrary deterministic continuous function defined on the real line and with values in an arbitrary fixed compact interval $[\underline{H},\overline{H}] \subset (1/\al,1)$;
\item $\Za{ds}$ is an independently scattered symmetric $\al$ stable ($\sas$) random measure on $\R$, with Lebesgue measure as its control measure. Many information on such random measures and the corresponding stochastic integrals can be found in \citep{SamTaq}.
\end{enumerate}
LMSM, denoted by $Y=\lbrace Y(t):t\in\R\rbrace$, of functional Hurst parameter $H(\cdot)$ is defined, for each $t\in\R$ as
\begin{equation}\label{defhaar:LMSM}
Y(t):=X(t,H(t)),
\end{equation}
where $X=\lbrace X(u,v):(u,v) \in \R\times (1/\al,1) \rbrace$, is the $\sas$ random field, such that for every $(u,v)\in\R \times (1/\al,1)$,
\begin{equation}
\label{defhaar:champstableX}
X(u,v):=\int_{\R}\Big\{ (u-s)_+^{v-1/\alpha} - (-s)_+^{v-1/\alpha} \Big\} \Za{ds},
\end{equation}
where for each real numbers $s$ and $\kappa$, 
\begin{equation}
\label{eqhaar:adpat+}
(s)_+^{\kappa}:=
\left\{
\begin{array}{l}
s^{\kappa}, \mbox{if $s\in (0,+\infty)$},\\
\\
0, \mbox{if $s\in (-\infty,0]$}.
\end{array}
\right.
\end{equation}
A modification of the high frequency part of $X$, is the $\sas$ stochastic field $\widetilde{X}_1=\{\widetilde{X}_1(u,v) : (u,v)\in \R \times (1/\alpha,1)\}$ defined for each $(u,v)\in\R\times (1/\al,1)$, as:
\begin{equation}
\label{eqhaar:hfX}
\widetilde{X}_1(u,v):=\int_{0}^{+\infty} (u-s)_+^{v-1/\alpha} \Za{ds}\,;
\end{equation}
a modification of the low frequency part of $X$, is the $\sas$ stochastic field $\widetilde{X}_2=\{\widetilde{X}_2(u,v) : (u,v)\in \R\times (1/\alpha,1)\}$ defined for each $(u,v)\in\R\times (1/\al,1)$, as:
\begin{equation}
\label{eqhaar:lfX}
\widetilde{X}_2(u,v):=\int_{-\infty}^{0} \Big\{ (u-s)_+^{v-1/\alpha} - (-s)_+^{v-1/\alpha} \Big\} \Za{ds}\,.
\end{equation}
It is worth noticing that the properties of these two fields are far from being completely similar. Also, observe that, in view of (\ref{defhaar:champstableX}), (\ref{eqhaar:hfX}) 
and (\ref{eqhaar:lfX}), one has for all $(u,v)\in \R \times (1/\alpha,1)$, almost surely,
$$
X(u,v)= \widetilde{X}_1(u,v)+\widetilde{X}_2(u,v).
$$

Note that for each $v \in (1/\al,1)$, the process $X(\cdot,v)=\lbrace X(u,v):u\in\R\rbrace$ is the LFSM of Hurst parameter $v$. Observe also that, in the particular case $\al=2$, LMSM reduces to the Multifractional Brownian Motion of functional Hurst parameter $H(\cdot)$.

In a certain way, this article is based on \citep{ayache2009series};  our main objective is to introduce a Haar expansion of the high and low frequency parts of the $\sas$ random field $X$, which generate LMSM. Then, to show that these series are convergent in a strong sense, namely, almost surely, in the space of continuous functions.  Finally, to determine their almost sure rates of convergence in the latter space. It is worth noticing that these representations of the high and low frequency parts of $X$, provide a new method for simulating the high and low frequency parts of LMSM. In the Gaussian case, many results concerning random wavelet series representations of such random models have been derived in literature (see e.g. \citep{Meyer1999,ayache2003rate,ayache2009series,lifshits2009haar}). 
In the Stable case, there are few results concerning wavelet series representations of LFSM or LMSM (see e.g. \citep{ayache2009linear,ayache2012lmsm}) 

In order to motivate this paper, let us recall the main result of \citep[Section 3]{ayache2012lmsm}. Let $\psi:\R\rightarrow\R$ be a 3 times continuously differentiable compactly supported Daubechies mother wavelet \citep{Dau92,Meyer90,Meyer92} and let $\mathcal{E}_{\gamma} (a,b,M):=\ce^{1}\big( [a,b], \ce^{\gamma}([-M,M],\R) \big)$ be the Banach space
of the Lipschitz functions defined on $[a,b]$ and with values in the H\"{o}lder space $\ce^{\gamma}([-M,M],\R)$, where $\ga$, $a$, $b$ and $M$ denote four arbitrary real numbers satisfying $M>0$, $1/\al<a<b<1$ and $0\le \ga <a-1/\al$. One has, for all $\o\in\Omega$, where $\Omega$ is a event of probability $1$,  in $\mathcal{E}_{\gamma} (a,b,M)$,
\begin{equation}
\label{eq1haar:descrip}
X(u,v,\o)=\lim_{n\rightarrow +\infty} \sum_{(j,k)\in D_{M,n}} 2^{-jv}\epsilon_{j,k}(\o)\big(\Psi (2^j u-k,v)-\Psi (-k,v)\big),
\end{equation}
where: 
\begin{itemize}
\item $D_{M,n}$ is the set of indices defined as,
\begin{equation}
\label{eq2haar:PWSEXN}
D_{M,n}:=\left\{(j,k)\in \Z^2 : |j| \leq n \mbox{ and } |k| \leq M 2^{n+1}\right\},
\end{equation} 
\item $\Psi$ is the smooth deterministic function defined for all $(x,v)\in \R\times (1/\alpha,1)$ as, 
\begin{equation}\label{daub:PSI}
\Psi(x,v):=\int_{\R} (x-s)_{+}^{v-1/\alpha} \psi(s) ds\,,
\end{equation}
\item $\lbrace \epsilon_{j,k} : (j,k)\in\Z^2 \rbrace$ is the sequence of the real-valued random variables defined as,
\begin{equation}
\label{daubdef:ejk}
\epsilon_{j,k} := 2^{j/\alpha}\int_{\R}  \psi(2^js-k) \Za{ds}\,;
\end{equation}
observe that, our assumption that $\Za{ds}$ is a $\sas$ random measure, implies that all the $\epsilon_{j,k}$'s have the same $\sas$ distribution of scale parameter equals $\|\psi\|_{L^{\al}(\R)}$.
\end{itemize}
Also recall that the following two results, are the two main ingredients of the proof of the fact that the convergence in (\ref{eq1haar:descrip}) holds in the space $\mathcal{E}_{\gamma} (a,b,M)$.
\begin{enumerate}
\item The sequence of real numbers $\lbrace \epsilon_{j,k}(\o) : (j,k)\in\Z^2 \rbrace$ satisfies (see \citep[Corollary 5]{ayache2009linear}), for every fixed arbitrarily small $\eta>0$,  
\begin{equation}
\label{eq1haar:omega0}
\big|\epsilon_{j,k}(\o)\big| \le C'(\o)\big(3+|j| \big)^{1/\alpha+\eta}\big(3+|k| \big)^{1/\alpha+\eta}, 
\end{equation}
where $C'$ is a positive and finite random variable only depending on $\eta$.
\item The function $\Psi$ as well as all its partial derivatives of any order, are well-localized in the variable $x$ uniformly in the variable $v\in [a,b]$, namely, for each $(p,q) \in \{ 0,1,2,3 \}\times\Z_+$, one has
\begin{equation}
\label{haar:localisation}
\sup_{(x,v)\in\R\times [a,b]} (3+|x|)^2 \big|(\partial_x^p \partial_v^q \Psi)(x,v)\big|< \infty.
\end{equation}
\end{enumerate}
The proofs of several results in \citep{ayache2012lmsm} testify that, the series representation in (\ref{eq1haar:descrip}) as well as its pathwise partial derivative with respect to $v$, are powerful tools for a fine study of path properties of the field $X=\{X(u,v):(u,v)\in \R\times (1/\alpha,1) \}$ and a corresponding LMSM defined above. However, this representation via Daubechies wavelets, has the following two drawbacks.
\begin{itemize}
\item The function $\Psi$ (see (\ref{daub:PSI})) and the random variables $\eps_{j,k}$ (see (\ref{daubdef:ejk})) cannot be defined by simple explicit formulas, since this is not the case for the Daubechies wavelet $\psi$ itself; therefore
(\ref{eq1haar:descrip}) can hardly provide an efficient simulation method of the field $X$ and $Y$ a corresponding LMSM.
\item In (\ref{eq1haar:descrip}), the high and the low frequency parts of $X$ (namely the fields $\widetilde{X}_1$ and $\widetilde{X}_2$ defined through (\ref{eqhaar:hfX}) and (\ref{eqhaar:lfX})) are not completely separated; basically, this comes from the fact that the diameter of the 
support of $\psi$, is strictly larger than $1$ (it is even much larger than $1$).
\end{itemize}

In order to avoid the latter two drawbacks, we replace $\psi$ by the Haar mother wavelet $h$ defined for all $s\in\R$, as:
\begin{equation}
\label{eqhaar:haarmwaveley}
h(s):=\one_{[0,1/2)}(s)-\one_{[1/2,1)}(s),
\end{equation}
where $\one_S$ is the indicator function of an arbitrary subset $S$ of $\R$. The continuously differentiable function $\th:\R\times (1/\al,1)\rightarrow\R$ is defined through (\ref{daub:PSI}) in which $\psi$ is replaced by $h$; it is worth noticing that, despite the fact that $\th$ will basically play the same role as $\Psi$, there is a considerable difference between both functions, indeed:
\begin{itemize}
\item on one hand, $\th$ has the advantage to be explicitly given by a simple formula, namely, in view of (\ref{eqhaar:haarmwaveley}), for all
$(x,v)\in\R\times (1/\al,1)$, one has,
\begin{align}
\label{haar:PSI}
\theta(x,v)& := \int_{\R} (x-s)_{+}^{v-1/\alpha} h(s) ds\\
 &=  \big(1+v-1/\al\big)^{-1} \Big\{ (x-1)_{+}^{1+v-1/\al}-2(x-1/2)_{+}^{1+v-1/\al}+ (x)_{+}^{1+v-1/\al} \Big\} \label{theta:00uv:alpha};
\end{align}
\item but, on the other hand, $\th$ is less regular than $\Psi$ and more importantly it fails to satisfy the "nice" localization property (\ref{haar:localisation}).
\end{itemize}

For each $(j,k)\in\Z^2$, we denote by $\lajk$ the $\sas$ random variable defined through (\ref{daubdef:ejk}) in which $\psi$ is replaced by $h$; in contrast with $\epsilon_{j,k}$, the random variable $\lajk$ is explicitly given by a simple formula, namely, in view of (\ref{eqhaar:haarmwaveley}), one has:
\begin{align}
\lajk &:=  2^{j/\al}\int_{\R} h(2^j s-k) \Za{\mathrm{ds}} \label{haardef:ejk}\\
& =  -2^{j/\alpha} \bigg(\Za{\frac{k}{2^j}}-2\Za{\frac{k+1/2}{2^j}}+\Za{\frac{k+1}{2^j}} \bigg),\label{haar:ejk}
\end{align}
where $\{\Za{t}:t\in\R\}$ is the $\sas$ L\'evy process with c\`adl\`ag paths which has been introduced at the very beginning of this introduction. Observe that the $\lajk$'s are identically distributed with a scale parameter equals $1$; also observe that for every fixed $j\in\Z$, $\{\lajk :k\in\Z\}$ is a sequence of independent 
random variables.

In the rest of this article, we always assume that $(u,v)$ belongs to the compact rectangle $[0,1]\times [a,b]$, where $a$ and $b$ are two arbitrary fixed real numbers satisfying $1/\al<a<b<1$;
typically one has $a=\min_{x\in [0,1]} H(x)$ and $b=\max_{x\in [0,1]} H(x)$ where $H(\cdot)$ is the continuous functional parameter of $Y$, the LMSM defined through (\ref{defhaar:LMSM}). The stochastic $\sas$ fields $X$, $X_1$ and $X_2$ are identified with their restrictions to $[0,1]\times [a,b]$.

Let us now introduce random series representations of $X_1$ and $X_2$ via Haar functions. On one hand, (\ref{eqhaar:hfX}), the fact that the sequence of functions:
$$
\big\{ \one_{[0,1]}(\cdot) \big\} \cup \left\{2^{j/2}h(2^j\cdot - k): j\in\Z_+ \mbox{ and } k\in\{0,\ldots, 2^j-1\} \right\},
$$
is an orthonormal basis of the Lebesgue Hilbert space $L^2\big([0,1]\big)$ (see \citep{Dau92,Meyer90,Meyer92}), standard computations, H\"older inequality, and a classical property
of the stochastic integral with respect to $\mathrm{Z}_{\alpha}$, imply that for all fixed $(u,v)\in [0,1]\times [a,b]$, the sequence of the $\sas$ random variables $\big(X_{1}^J(u,v)\big)_{J\in\N}$,
defined as:
\begin{equation}
\label{eqhaar:partialsum-hf}
X_{1}^J (u,v)=\frac{u^{1+v-1/\al}}{1+v-1/\al} \Za{1} +\sum_{j=0}^{J-1} 2^{-jv}\sum_{k=0}^{2^j -1}\lajk\th (2^j u-k,v),
\end{equation}
converges in probability to the random variable $X_1(u,v)$, when $J$ goes to $+\infty$.
On the other hand, (\ref{eqhaar:lfX}), the fact that the sequence of functions:
$$
\left\{2^{j/\al}h(2^j\cdot + k): (j,k)\in\Z\times\N\right\},
$$
is an unconditional basis of the Lebesgue space $L^{\al}\big((-\infty,0]\big)$ (see \citep{Dau92,Meyer90,Meyer92}), standard computations, and a classical property
of the stochastic integral with respect to $\mathrm{Z}_{\alpha}$, imply that for all fixed $(u,v)\in [0,1]\times [a,b]$, the sequence of the $\sas$ random variables $\big(X_{2}^J(u,v)\big)_{J\in\N}$,
defined as:
\begin{equation}
\label{eqhaar:partialsum-lf}
X_{2}^J (u,v)=\sum_{j=1-J}^{J-1} 2^{-jv}\sum_{k=1}^{2^{J-|j|}} \lajnk\big(\th (2^j u+k,v)-\th (k,v)\big),
\end{equation}
converges in probability to the random variable $X_2(u,v)$, when $J$ goes to $+\infty$.
It is clear that for every $J\in\N$, the paths of the $\sas$ fields $X_{1}^J =\{X_{1}^J (u,v):(u,v)\in [0,1]\times [a,b]\}$ and $X_{2}^J =\{X_{2}^J (u,v):(u,v)\in [0,1]\times [a,b]\}$ belong to $\ce:=\ce([0,1]\times [a,b],\R)$, the Banach space of the real-valued continuous functions on the rectangle $[0,1]\times [a,b]$ equipped 
with the usual supremum norm, denoted by $\|\cdot\|_{\ce}$. A natural question one can address is that, whether or not, the sequences of continuous random functions $(X_{1}^J)_{J\in\N}$ and 
$(X_{2}^J)_{J\in\N}$, almost surely converge in the space $\ce$; assume for a while that the answer to the question is positive and denote by $X_1$ and $X_2$ the limits
of these two sequences, then the $\sas$ fields $X_{1} =\{X_{1} (u,v):(u,v)\in [0,1]\times [a,b]\}$ and $X_{2} =\{X_{2} (u,v):(u,v)\in [0,1]\times [a,b]\}$ are modifications with almost surely continuous paths, respectively of the high and the low frequency parts of the field $X =\{X (u,v):(u,v)\in [0,1]\times [a,b]\}$. 

The main difficulty to show that the sequences $(X_{1}^J)_{J\in\N}$ and $(X_{2}^J)_{J\in\N}$ are almost surely convergent in the space $\ce$ is that the function $\th$ (see (\ref{haar:PSI})) is a badly localized function in the variable $x$; actually when $v\in [a,b]$ is fixed and $x$ goes to $+\infty$, then $\th (x,v)$ goes to $0$ enough fast, it vanishes at the same rate as $x^{v-1/\al-1}$. Abel transforms is the important tool, which allow us to overcome this difficulty.

The two main results of this paper are the following.

\begin{thm}
\label{thhaar:thmain1}
Let $\Omega_{0}^{**}$ be the event of probability 1 introduced in Proposition~\ref{prophaar:estimljk}. Then, for each $\o\in\Omega_{0}^{**}$, the sequence of the continuous functions
$(X_{1}^J (\cdot,\cdot,\o))_{J\in\N}$ defined through (\ref{eqhaar:partialsum-hf}), converges in the space $\ce:=\ce([0,1]\times [a,b],\R)$ to a limit denoted by 
$X_1(\cdot,\cdot,\o)$; moreover, one has for all fixed $\eta>0$, 
\begin{equation}
\label{thmain1:eq1haar}
\big\| X_1(\cdot,\cdot,\o)-X_{1}^J (\cdot,\cdot,\o)\big\|_{\ce} =\mathcal{O}\left(2^{-J(a-1/\al)} J^{2/\al+\eta}\right).
\end{equation} 
Notice that the $\sas$ field $X_{1} =\{X_{1} (u,v):(u,v)\in [0,1]\times [a,b]\}$ is a modification with almost surely continuous paths, of the high frequency part of the $\sas$ field $X =\{X (u,v):(u,v)\in [0,1]\times [a,b]\}$ which generates LMSM's. 
\end{thm}

\begin{thm}
\label{thhaar:thmain2}
Let $\Omega_{1}^{**}$ and $\Omega_{2}^{**}$ be the events of probability 1 introduced in Propositions~\ref{prophaar:estimljk2}~and~\ref{prophaar:estimljk2-}. Then, for each $\o\in\Omega_{1}^{**}\cap \Omega_{2}^{**}$ (notice that the event $\Omega_{1}^{**}\cap \Omega_{2}^{**}$ is of probability $1$), the sequence of the continuous functions
$(X_{2}^J (\cdot,\cdot,\o))_{J\in\N}$ defined through (\ref{eqhaar:partialsum-lf}), converges in the space $\ce:=\ce([0,1]\times [a,b],\R)$ to a limit denoted by 
$X_2 (\cdot,\cdot,\o)$; moreover, one has for all fixed $\eta>0$, 
\begin{equation}
\label{thmain2:eq1haar}
\big\| X_2(\cdot,\cdot,\o)-X_{2}^J (\cdot,\cdot,\o)\big\|_{\ce} =\mathcal{O}\left(2^{-J(1-b)} J^{1/\al+\eta}\right).
\end{equation} 
Notice that the $\sas$ field $X_{2} =\{X_{2} (u,v):(u,v)\in [0,1]\times [a,b]\}$ is a modification with almost surely continuous paths, of the low frequency part of the $\sas$ field $X =\{X (u,v):(u,v)\in [0,1]\times [a,b]\}$ which generates LMSM's.
\end{thm}

Notice, in Theorem \ref{thhaar:thmain1}, that the speed of convergence of the series $(X_1^J)_{J\in \N}$ is better when the parameter $v$ takes values close to $1$; conversely, in Theorem \ref{thhaar:thmain2}, the rate of convergence of the series $(X_2^J)_{J\in \N}$ is better when $v$ takes values close to $1/\al$.

The remaining of the paper is structured in the following way: Section \ref{sechaar:thmain1} is devoted to the proof of Theorem \ref{thhaar:thmain1} and Section \ref{sechaar:thmain2} to that of Theorem \ref{thhaar:thmain2}. Finally, some simulations of LMSM are given in Section \ref{sechaar:simulation}.

\section{Proof of Theorem~\ref{thhaar:thmain1}}\label{sechaar:thmain1}

Let $\Th:\R\times (1/\al,1)\rightarrow\R$ be the continuously differentiable function defined for every $(x,v)\in\R\times (1/\al,1)$ as:
\begin{equation}
\label{eq2haar:abelhf}
\Th (x,v):=\th (x,v)-\th(x-1,v).
\end{equation}
For all $(j,k)\in\Z_{+}^2$, let $\ljk$ be the $\sas$ random variable of scale parameter $(k+1)^{1/\al}$ defined as:
\begin{equation}
\label{eq2haarbis:abelhf}
\ljk:=\sum_{m=0}^{k} \zeta_{j,m}.
\end{equation}

\begin{lemma}
\label{lemhaar:abelhf}
One has for each $j\in\Z_+$ and $(u,v)\in [0,1]\times [a,b]$,
\begin{equation}
\label{eq3haar:abelhf}
\sum_{k=0}^{2^j -1}\lajk\th (2^j u-k,v)=\lambda_{j,2^j-1}\th(2^j u-2^j+1,v)+\sum_{k=0}^{2^j -2}\ljk\Th(2^j u-k,v),
\end{equation}
with the convention that $\sum_{k=0}^{-1}\lambda_{0,k}\Th(u-k,v)=0$.
\end{lemma}

\begin{proof}[Proof of Lemma~\ref{lemhaar:abelhf}] It is clear that (\ref{eq3haar:abelhf}) holds when $j=0$, so let us assume that $j\ge 1$. Using
(\ref{eq2haarbis:abelhf}) and (\ref{eq2haar:abelhf}), one obtains that,
\begin{eqnarray*}
\sum_{k=0}^{2^j -1}\lajk\th (2^j u-k,v)&=& \lambda_{j,0}\th(2^j u,v)+\sum_{k=1}^{2^j-1} (\ljk-\lambda_{j,k-1})\th (2^j u-k,v)\\
&=& \sum_{k=0}^{2^j -1} \ljk\th (2^j u-k,v)-\sum_{k=0}^{2^j-2}\ljk\th (2^j u-k-1,v)\\
&=& \lambda_{j,2^j-1}\th(2^j u-2^j+1,v)+\sum_{k=0}^{2^j -2}\ljk\Th(2^j u-k,v).
\end{eqnarray*}
\end{proof}

Let us now provide a rather sharp estimate of the asymptotic behavior of the sequence of random variables $\{\ljk:(j,k)\in\Z_{+}^2\}$.
\begin{Prop}
\label{prophaar:estimljk}
There exists an event of probability 1, denoted by $\Omega_0^{**}$, such that for every fixed real number $\eta >0$, one has, for all $\omega \in \Omega_0^{**}$ and for each $(j,k)\in\Z_{+}^2$,
\begin{equation}
\label{eq1:haar:estimljk}
\big| \ljk (\o)\big| \leq C(\omega) \big(1+j \big)^{1/\alpha}\log^{1/\alpha+\eta}\big(3+ j\big) \, \big(1+k \big)^{1/\alpha} \log^{1/\alpha+\eta}\big(3+ k\big), 
\end{equation}
where $C$ is a positive and finite random variable only depending on $\eta$.
\end{Prop}
In order to prove Proposition~\ref{prophaar:estimljk}, we need two preliminary results.
\begin{lemma}
\label{lem:distljk}
For each fixed $j\in\Z_+$, the $\sas$ process $\{\ljk: k\in\Z_{+}\}$ has the same finite dimensional distributions as the process $\{\Za{k+1}: k\in\Z_+\}$; recall 
that $\{\Za{t}: t\in\R_+\}$ is a $\sas$ L\'evy process with c\`adl\`ag paths.
\end{lemma}

\begin{proof}[Proof of Lemma~\ref{lem:distljk}] Let $\{\de_m: m\in\Z_+\}$ be the sequence of the independent and identically distributed $\sas$ random variables with a scale parameter equals $1$, defined for all $m\in\Z_+$ as:
\begin{equation}
\label{eq2quat:haar:estimljk}
\de_m=\Za{m+1}-\Za{m}.
\end{equation}
It follows from (\ref{eq2quat:haar:estimljk}) and the equality $Z_{\al}(0)=0$, that one has for each $k\in\Z_+$,
\begin{equation}
\label{eq3quat:haar:estimljk}
\Za{k+1}=\sum_{m=0}^k \de_m.
\end{equation}
Then, combining (\ref{eq3quat:haar:estimljk}) with (\ref{eq2haarbis:abelhf}), and using the fact that for each fixed $j\in\Z_+$, $\{\zeta_{j,m}:m\in\Z_+\}$ is a sequence of  independent and identically distributed $\sas$ random variables with a scale parameter equals $1$, one gets the lemma.
\end{proof}

\begin{lemma}
\label{lem:haarkin}
Let $\eta>0$ be arbitrary and fixed. We set,
\begin{equation}
\label{eq1:haarkin}
M_{\eta,\al}^*:=\sup\left\{\frac{|\Za{t}|}{t^{1/\al}\log^{1/\al+\eta}(2+t)}:t\in [1,+\infty)\right\}.
\end{equation}
Then $M_{\eta,\al}^*$ is an almost surely finite random variable; moreover there is a constant $c>0$ such that for all real number $y\ge 1$, one has,
\begin{equation}
\label{eq2:haarkin}
\PR\left(M_{\eta,\al}^*>y\right)\le c y^{-\al}.
\end{equation}
\end{lemma}

\begin{proof}[Proof of Lemma~\ref{lem:haarkin}] The fact that $M_{\eta,\al}^*$ is an almost surely finite random variable, has been derived in \citep{Khin38}. The inequality (\ref{eq2:haarkin}) can be obtained by applying \citep[Theorem~10.5.1]{SamTaq}, to the almost surely  bounded $\sas$ processes 
$
\left\{\frac{\Za{t}}{t^{1/\al}\log^{1/\al+\eta}(2+t)}:t\in [1,+\infty)\right\}
$
and 
$
\left\{-\frac{\Za{t}}{t^{1/\al}\log^{1/\al+\eta}(2+t)}:t\in [1,+\infty)\right\}.
$
\end{proof} 

Now, we are in position to prove Proposition~\ref{prophaar:estimljk}.

\begin{proof}[Proof of Proposition~\ref{prophaar:estimljk}] For each fixed $j\in\Z_+$, let $\mu_{\eta,\al}^{j,*}$ be the random variable defined as:
\begin{equation}
\label{eq2:haar:estimljk}
\mu_{\eta,\al}^{j,*}:=\sup\left\{\frac{|\ljk|}{(1+k)^{1/\al}\log^{1/\al+\eta}(3+k)}: k\in\Z_+\right\}.
\end{equation}
In view of Lemma~\ref{lem:distljk}, one has, for all $j\in\Z_+$, 
\begin{equation}
\label{eq2bis:haar:estimljk}
\mu_{\eta,\al}^{j,*}\stackrel{(d)}{=}\nu,
\end{equation}
where $\stackrel{(d)}{=}$ means equality in distribution, and $\nu$ is the random variable, defined as:
\begin{equation}
\label{eq3:haar:estimljk}
\nu:=\sup\left\{\frac{|\Za{k+1}|}{(1+k)^{1/\al}\log^{1/\al+\eta}(3+k)}: k\in\Z_+\right\}.
\end{equation}
Notice that Lemma~\ref{lem:haarkin} implies that $\nu$ is almost surely finite; moreover, thanks to (\ref{eq2bis:haar:estimljk}), the latter property is also satisfied by
$\mu_{\eta,\al}^{j,*}$, for any arbitrary $j\in\Z_+$. Next, using (\ref{eq2bis:haar:estimljk}), (\ref{eq3:haar:estimljk}), (\ref{eq1:haarkin}) and (\ref{eq2:haarkin}), one gets that,
\begin{align*}
& \sum_{j=0}^{+\infty} \PR\left(\mu_{\eta,\al}^{j,*}> \big(1+j \big)^{1/\alpha}\log^{1/\alpha+\eta}\big(3+ j\big)\right) \\
&= \sum_{j=0}^{+\infty} \PR\left(\nu> \big(1+j \big)^{1/\alpha}\log^{1/\alpha+\eta}\big(3+ j\big)\right)\\
& \le  \sum_{j=0}^{+\infty}\PR\left(M_{\eta,\al}^*> \big(1+j \big)^{1/\alpha}\log^{1/\alpha+\eta}\big(3+ j\big)\right)\\
&\le  c\sum_{j=0}^{+\infty}\big(1+j \big)^{-1}\log^{-1-\al\eta}\big(3+ j\big)<\infty;
\end{align*}
thus the proposition results from Borel-Cantelli Lemma as well as from the fact that $\mu_{\eta,\al}^{j,*}$ is almost surely finite for each $j\in\Z_+$.
\end{proof}

The following proposition provides sharp estimates of the rate of vanishing of $\th(x,v)$ and $\Th(x,v)$ (see (\ref{theta:00uv:alpha}) and (\ref{eq2haar:abelhf})) when $x$ goes to infinity.
\begin{Prop}
\label{prophaar:localTh}
\begin{itemize}
\item[(i)] For each $(x,v)\in (-\infty,0]\times (1/\al,1)$, one has, 
$$\th(x,v)=\Th(x,v)=0.$$
\item[(ii)] There exists a constant $c>0$, such that for all $(x,v)\in (0,+\infty)\times [a,b]$, one has,
$$
|\th(x,v)|\le c(1+x)^{v-1/\al-1}.
$$
\item[(iii)] There is a constant $c'>0$, such that for all $(x,v)\in (0,+\infty)\times [a,b]$, one has,
$$
|\Th(x,v)|\le c'(1+x)^{v-1/\al-2}.
$$
\end{itemize} 
\end{Prop}

\begin{proof}[Proof of Proposition~\ref{prophaar:localTh}] In view of (\ref{theta:00uv:alpha}), (\ref{eq2haar:abelhf}) and (\ref{eqhaar:adpat+}), it is clear that 
Part $(i)$ of the proposition holds. Let us prove the two other parts of it. Observe that the fact that $(x,v)\mapsto (1+x)^{1-v-1/\al}\,\th(x,v)$ and 
$(x,v)\mapsto (1+x)^{2-v-1/\al}\,\Th(x,v)$ are continuous functions on the compact rectangle $[0,4]\times [a,b]$, implies that,
\begin{equation}
\label{eq1haar:localTh}
c_1:=\sup\left\{(1+x)^{1-v-1/\al}|\th(x,v)|:(x,v)\in [0,4]\times [a,b]\right\}<\infty 
\end{equation}
and
\begin{equation}
\label{eq2haar:localTh}
c_2:=\sup\left\{(1+x)^{2-v-1/\al}|\Th(x,v)|:(x,v)\in [0,4]\times [a,b]\right\}<\infty.
\end{equation}
From now on, we assume that $(x,v)\in (4,+\infty)\times [a,b]$. Observe that, in view of (\ref{theta:00uv:alpha}) and (\ref{eqhaar:adpat+}), one has,
\begin{equation}
\label{eq3haar:localTh}
 \th(x,v)=\frac{x^{1+v-1/\al}}{1+v-1/\al}\left\{\big(1-x^{-1}\big)^{1+v-1/\al}-2\big(1-(2x)^{-1}\big)^{1+v-1/\al}+ 1 \right\}.
 \end{equation}
Next, let us show that there are two constants $c_3>0$ and $c_4>0$, such that for all $(z,v)\in [0,2^{-1}]\times [a,b]$, one has,
\begin{equation}
\label{eq4haar:localTh}
\left|(1-z)^{1+v-1/\al}-1+(v-1/\al+1)z\right|\le c_3 z^2
\end{equation}
and 
\begin{equation}
\label{eq5haar:localTh}
\left|(1-z)^{1+v-1/\al}-1+(v-1/\al+1)z-\frac{(v-1/\al+1)(v-1/\al)}{2}z^2\right|\le c_4 z^3.
\end{equation}
Observe that (\ref{eq4haar:localTh}) easily results from (\ref{eq5haar:localTh}), so we only need to prove that the latter inequality holds. Applying, for each fixed $v\in [a,b]$, Taylor-Lagrange formula, to the function $y\mapsto (1-y)^{1+v-1/\al}$, on the interval $[0,z]$, one gets that 
\begin{eqnarray}
\label{eq6haar:localTh}
&& (1-z)^{1+v-1/\al}-1+(v-1/\al+1)z-\frac{(v-1/\al+1)(v-1/\al)}{2}z^2\nonumber\\
&& =-\frac{(v-1/\al+1)(v-1/\al)(v-1/\al-1)}{6}(1-\xi)^{v-1/\al-2}\, z^3,
\end{eqnarray}
where $\xi\in [0,z]\subset [0,2^{-1}]$, then (\ref{eq5haar:localTh}) easily follows from (\ref{eq6haar:localTh}). Next, using the triangle inequality and (\ref{eq4haar:localTh}) (in the case where $z=x^{-1}$ and also in the case where $z=(2x)^{-1}$), one gets that,
\begin{eqnarray}
\label{eq7haar:localTh}
&& \left|\big(1-x^{-1}\big)^{1+v-1/\al}-2\big(1-(2x)^{-1}\big)^{1+v-1/\al}+ 1 \right|\nonumber\\
&& =\left|\big(1-x^{-1}\big)^{1+v-1/\al}-1+(v-1/\al+1)x^{-1}\right.\nonumber\\
&&  \hspace{3.5cm}\left. -2\big(1-(2x)^{-1}\big)^{1+v-1/\al}+ 2-2(v-1/\al+1)(2x)^{-1}\right|\nonumber\\
&& \le \left|\big(1-x^{-1}\big)^{1+v-1/\al}-1+(v-1/\al+1)x^{-1}\right|\nonumber\\
&& \hspace{3.5cm}+2\left|\big(1-(2x)^{-1}\big)^{1+v-1/\al}- 1+(v-1/\al+1)(2x)^{-1}\right|\nonumber\\
&& \le  c_5 x^{-2},
\end{eqnarray}
where $c_5:=(3/2)c_3$. Next, putting together (\ref{eq1haar:localTh}), (\ref{eq3haar:localTh}) and (\ref{eq7haar:localTh}), one obtains Part~$(ii)$ of the proposition. Let us now prove that 
Part~$(iii)$ of it, holds; to this end, we set 
\begin{equation}
\label{eqhaar:defdl}
\mbox{$d_0=1$, $d_1=-2$, $d_2=0$, $d_3=2$ and $d_4=-1$.} 
\end{equation}
Standard computations allow to show that, for all $m\in\{0,1,2\}$, 
\begin{equation}
\label{eq8haar:localTh}
\sum_{l=0}^{4} l^m d_l=0,
\end{equation}
with the convention that $0^0:=1$; moreover, in view of (\ref{eq2haar:abelhf}), (\ref{theta:00uv:alpha}) and (\ref{eqhaar:adpat+}), for each $(x,v)\in (4,+\infty)\times [a,b]$, one has,
\begin{eqnarray}
\label{eq9haar:localTh}
\Th(x,v)&=& \big(1+v-1/\al\big)^{-1}\left\{\sum_{l=0}^{4} d_l\big(x-l/2\big)^{1+v-1/\al}\right\}\nonumber\\
&=& \frac{x^{1+v-1/\al}}{\big(1+v-1/\al\big)}\left\{\sum_{l=0}^{4} d_l \Big (1-\frac{l}{2x}\Big)^{1+v-1/\al}\right\}.
\end{eqnarray} 
Next, using (\ref{eq8haar:localTh}), the triangle inequality, and (\ref{eq5haar:localTh}) in which one takes $z=l/(2x)$, it follows that,
\begin{eqnarray}
\label{eq10haar:localTh}
&& \left|\sum_{l=0}^{4} d_l \Big (1-\frac{l}{2x}\Big)^{1+v-1/\al}\right|\nonumber\\
&& = \left|\sum_{l=0}^{4} d_l\left\{\Big(1-\frac{l}{2x}\Big)^{1+v-1/\al}\right.\right. \nonumber\\
&& \hspace{2.5cm} \left.\left.-1+(v-1/\al+1)\frac{l}{2x}-\frac{(v-1/\al+1)(v-1/\al)}{2}\big(\frac{l}{2x}\big)^2\right\}\right|\nonumber\\
&& \le c_6 x^{-3},
\end{eqnarray}
where $c_6:=8^{-1}c_4 \big(\sum_{l=1}^4 l^3 |d_l|\big)$. Finally, putting together (\ref{eq2haar:localTh}), (\ref{eq9haar:localTh}) and (\ref{eq10haar:localTh}), one 
obtains Part $(iii)$ of the proposition.
\end{proof}

\begin{lemma}
\label{lemhaar:supTh}
One has 
\begin{equation}
\label{eq1haar:supTh}
M:=\sup\left\{\sum_{k\in\Z} \big|\Th (x-k, v)\big|: (x,v)\in \R\times [a,b]\right\}<\infty.
\end{equation}
\end{lemma}

\begin{proof}[Proof of Lemma~\ref{lemhaar:supTh}] Observe that, for each fixed $v\in [a,b]$, the function $$x\mapsto \sum_{k\in\Z} \big|\Th (x-k, v)\big|,$$ 
defined on $\R$ and a priori taking its values in $\R\cup\{+\infty\}$, is $1$-periodic; therefore it is sufficicent to show that (\ref{eq1haar:supTh}) holds 
when $(x,v)\in [-1/2,1/2]\times [a,b]$. Using Parts $(i)$ and $(iii)$ of Proposition~\ref{prophaar:localTh}, as well as, the triangle inequality, one gets that, 
\begin{eqnarray*}
\sum_{k\in\Z} \big|\Th (x-k, v)\big|&\le & c' \sum_{k\in\Z} \big(1+|x-k|\big)^{b-1/\al-2}\\
&\le & c'\sum_{k\in\Z} \big(1+|k|-|x|\big)^{b-1/\al-2}\\
&\le & c'\sum_{k\in\Z} \big(2^{-1}+|k|\big)^{b-1/\al-2}<\infty.
\end{eqnarray*}
\end{proof}

Now, we are in position to show that Theorem~\ref{thhaar:thmain1} holds.

\begin{proof}[Proof of Theorem~\ref{thhaar:thmain1}] Let $\o\in\Omega_{0}^{**}$ be arbitrary and fixed; recall that $\Omega_{0}^{**}$ is the event of probability 1 introduced in Proposition~\ref{prophaar:estimljk}. Let us first show that the sequence of the continuous functions
$(X_{1}^J (\cdot,\cdot,\o))_{J\in\N}$ defined through (\ref{eqhaar:partialsum-hf}), is a Cauchy sequence in $\ce:=\ce([0,1]\times [a,b],\R)$, the space of the real-valued continuous functions over $[0,1]\times [a,b]$, equipped with the usual supremum norm, denoted by $\|\cdot\|_{\ce}$. Let $\eta>0$ be arbitrary and fixed, using, (\ref{eqhaar:partialsum-hf}), the triangle 
inequality, (\ref{eq3haar:abelhf}), (\ref{eq1:haar:estimljk}), Parts $(i)$ and $(ii)$ of Proposition~\ref{prophaar:localTh}, and (\ref{eq1haar:supTh}), one gets that,
for all $(J,Q)\in\N^2$, 
\begin{align}
\label{thmain1:eq2haar}
& \big\|X_{1}^{J+Q} (\cdot,\cdot,\o)- X_{1}^{J}(\cdot,\cdot,\o)\big\|_{\ce} =\Big\|\sum_{j=J}^{J+Q-1} 2^{-j\cdot}\sum_{k=0}^{2^j -1}\lajk\th (2^j \cdot-k,\cdot)\Big\|_{\ce}\nonumber\\
& \le \sum_{j=J}^{J+Q-1} 2^{-ja}\Big\|\sum_{k=0}^{2^j -1}\lajk\th (2^j \cdot-k,\cdot)\Big\|_{\ce}\nonumber\\
& = \sum_{j=J}^{J+Q-1} 2^{-ja}\Big\|\lambda_{j,2^j-1}\th(2^j \cdot-2^j+1,\cdot)+\sum_{k=0}^{2^j -2}\ljk\Th(2^j \cdot-k,\cdot)\Big\|_{\ce}\nonumber\\
& \le \sum_{j=J}^{J+Q-1} 2^{-ja}\left(|\lambda_{j,2^j-1}|\big\|\th(2^j \cdot-2^j+1,\cdot)\big\|_{\ce}+ \right. \nonumber \\
&  \hspace{5cm} \left. \left(\max_{0\le l\le 2^j -2} |\ljk|\right)
\Big\|\sum_{k=0}^{2^j -2}|\Th(2^j \cdot-k,\cdot)|\Big\|_{\ce}\right)\nonumber\\
& \le C(\o)(c+M)\sum_{j=J}^{J+Q-1} 2^{-j(a-1/\al)}(2+j)^{2/\al+\eta}\log^{1/\al+\eta}(3+j)\nonumber\\
& \le C(\o)(c+M)\sum_{j=J}^{+\infty} 2^{-j(a-1/\al)}(2+j)^{2/\al+\eta}\log^{1/\al+\eta}(3+j).
\end{align}
It follows from (\ref{thmain1:eq2haar}) that $(X_{1}^J (\cdot,\cdot,\o))_{J\in\N}$ is a Cauchy sequence in $\ce$, and that $X_{1}(\cdot,\cdot,\o)$, its limit, satisfies 
for all $J\in\N$,
\begin{equation}
\label{thmain1:eq3haar}
\big\|X_{1}(\cdot,\cdot,\o)- X_{1}^{J}(\cdot,\cdot,\o)\big\|_{\ce}\le C(\o)(c+M)\sum_{j=J}^{+\infty} 2^{-j(a-1/\al)}(2+j)^{2/\al+\eta}\log^{1/\al+\eta}(3+j).
\end{equation}
Next, let us show that there exists a constant $c_1>0$, such that for all $J\in\N$, one has,
\begin{equation}
\label{thmain1:eq4haar}
\sum_{j=J}^{+\infty} 2^{-j(a-1/\al)}(2+j)^{2/\al+\eta}\log^{1/\al+\eta}(3+j)\le c_{1} 2^{-J(a-1/\al)}(2+J)^{2/\al+\eta}\log^{1/\al+\eta}(3+J).
\end{equation}
This is the case since,
\begin{align*}
& \sum_{j=J}^{+\infty} 2^{-j(a-1/\al)}(2+j)^{2/\al+\eta}\log^{1/\al+\eta}(3+j)\\
& = \sum_{j=0}^{+\infty} 2^{-(j+J)(a-1/\al)}(2+j+J)^{2/\al+\eta}\log^{1/\al+\eta}(3+j+J)\\
& \le \sum_{j=0}^{+\infty} 2^{-(j+J)(a-1/\al)}(2+j)^{2/\al+\eta}(2+J)^{2/\al+\eta}\log^{1/\al+\eta}\big\{(3+j)(3+J)\big\}\\
& \le 2^{-J(a-1/\al)}(2+J)^{2/\al+\eta}\sum_{j=0}^{+\infty}2^{-j(a-1/\al)}(2+j)^{2/\al+\eta}\Big(\log(3+j)+\log(3+J)\Big)^{1/\al+\eta}\\
& \le c_{1} 2^{-J(a-1/\al)}(2+J)^{2/\al+\eta}\log^{1/\al+\eta}(3+J),
\end{align*}
where 
$$
c_1:=\sum_{j=0}^{+\infty}2^{-j(a-1/\al)}(2+j)^{2/\al+\eta}\Big(\log(3+j)+1\Big)^{1/\al+\eta}.
$$
Finally, combining (\ref{thmain1:eq3haar}) with (\ref{thmain1:eq4haar}), one gets (\ref{thmain1:eq1haar}) in which $\eta$ is replaced by $2\eta$.
\end{proof}

\section{Proof of Theorem~\ref{thhaar:thmain2}}
\label{sechaar:thmain2}

\subsection{Study of the part $j\ge 0$ of the series}
\label{subsechaar:jpositive}

\begin{thm}
\label{th1haar:thmain2}
For each $J\in\N$, we denote by $X_{2,+}^J =\{X_{2,+}^J(u,v): (u,v)\in[0,1]\times [a,b]\}$ the $\sas$ field with paths in $\ce$, defined, for every 
$(u,v)\in[0,1]\times [a,b]$, as:
\begin{equation}
\label{eqhaar:partialsum-lf+}
X_{2,+}^J (u,v)=\sum_{j=0}^{J-1} 2^{-jv}\sum_{k=1}^{2^{J-j}} \lajnk\big(\th (2^j u+k,v)-\th (k,v)\big).
\end{equation}
Let $\Omega^{**}_{1}$ be the event of probability 1 introduced in Proposition~\ref{prophaar:estimljk2}. Then for all $\o\in\Omega^{**}_{1}$,  $\big(X_{2,+}^J(\cdot,\cdot,\o)\big)_{J\in\N}$ is a Cauchy sequence in $\ce$, moreover, its limit $X_{2,+}(\cdot,\cdot,\o)$, satisfies, for 
each fixed $\eta>0$,
\begin{equation}
\label{eq1:thmain2+}
\big\|X_{2,+}(\cdot,\cdot,\o)-X_{2,+}^J(\cdot,\cdot,\o)\big\|_{\ce}=\mathcal{O} \left (2^{-J(1-b)}J^{1/\al+\eta}\right)
\end{equation}
\end{thm}

In order to show that $\big(X_{2,+}^J(\cdot,\cdot,\o)\big)_{J\in\N}$ is a Cauchy sequence, one needs to appropriately bound the quantity $\big\|X_{2,+}^{J+Q}(\cdot,\cdot,\o)-X_{2,+}^{J}(\cdot,\cdot,\o)\big\|_{\ce}$, for all $(J,Q)\in\N^2$. Observe that, in view of (\ref{eqhaar:partialsum-lf+}), one has for every 
$(u,v)\in[0,1]\times [a,b]$,
\begin{equation}
\label{eqhaar:decompX2+}
X_{2,+}^{J+Q}(u,v,\o)-X_{2,+}^{J}(u,v,\o)=A_{2,+}^{J,Q}(u,v,\o)+B_{2,+}^{J,Q}(u,v,\o),
\end{equation}
where:
\begin{equation}
\label{eqhaar:defA2+}
A_{2,+}^{J,Q}(u,v,\o):=\sum_{j=0}^{J-1} 2^{-jv}\sum_{k=2^{J-j}+1}^{2^{J+Q-j}} \lajnk(\o)\big(\th (2^j u+k,v)-\th (k,v)\big),
\end{equation}
and 
\begin{equation}
\label{eqhaar:defB2+}
B_{2,+}^{J,Q}(u,v,\o):=\sum_{j=J}^{J+Q-1} 2^{-jv}\sum_{k=1}^{2^{J+Q-j}} \lajnk(\o)\big(\th (2^j u+k,v)-\th (k,v)\big).
\end{equation}
For all $(j,k)\in\Z_{+}\times\N$, let $\ljnk$ be the $\sas$ random variable of scale parameter $k^{1/\al}$ defined as:
\begin{equation}
\label{eq2haar:abellf+}
\ljnk:=\sum_{m=1}^{k} \zeta_{j,-m}.
\end{equation}

The proof of the following lemma is similar to that of Lemma~\ref{lemhaar:abelhf}.

\begin{lemma}
\label{lemhaar:abellf+}
Recall that the function $\Th$ has been introduced in (\ref{eq2haar:abelhf}). Let $(J,Q)\in\N^2$ and $(u,v)\in [0,1]\times [a,b]$ be arbitrary and fixed.
\begin{itemize}
\item[(i)] For each $j\in\{0,\ldots,J-1\}$, one has,
\begin{align}
\label{eq3haar:abellf+}
& \sum_{k=2^{J-j}+1}^{2^{J+Q-j}} \lajnk\th (2^j u+k,v)\\
&=\lambda_{j,-2^{J+Q-j}}\th (2^j u+2^{J+Q-j},v)-\lambda_{j,-2^{J-j}}\th (2^j u+2^{J-j}+1,v) \nonumber \\
& -\sum_{k=2^{J-j}+2}^{2^{J+Q-j}}\lambda_{j,-(k-1)}\Th(2^j u+k,v).\nonumber
\end{align}
\item[(ii)] For each $j\in\{J,\ldots,J+Q-1\}$, one has
\begin{eqnarray}
\label{eq4haar:abellf+}
&& \sum_{k=1}^{2^{J+Q-j}} \lajnk\th (2^j u+k,v)\\
&& =\lambda_{j,-2^{J+Q-j}}\th (2^j u+2^{J+Q-j},v)-\sum_{k=2}^{2^{J+Q-j}}\lambda_{j,-(k-1)}\Th(2^j u+k,v).\nonumber
\end{eqnarray}
\end{itemize}
\end{lemma}

The proof of the following proposition is similar to that of Proposition~\ref{prophaar:estimljk}.

\begin{Prop}
\label{prophaar:estimljk2}
There exists an event of probability 1, denoted by $\Omega_1^{**}$, such that for every fixed real number $\eta >0$, one has, for all $\omega \in \Omega_1^{**}$ and for each $(j,k)\in\Z_{+}\times\N$,
\begin{equation}
\label{eq1:haar:estimljk2}
\big| \ljnk (\o)\big| \leq C(\omega) \big(1+j \big)^{1/\alpha}\log^{1/\alpha+\eta}\big(3+ j\big)\,  k^{1/\alpha} \log^{1/\alpha+\eta}\big(2+ k\big), 
\end{equation}
where $C$ is a positive and finite random variable only depending on $\eta$.
\end{Prop}

The following propostion easily results from Proposition~\ref{prophaar:localTh}.
\begin{Prop}
\label{prophaar:localTh2}
There exist two constants $c>0$ and $c'>0$, such that for all $(j,k)\in\Z_{+}\times\N$ and $(u,v)\in [0,1]\times [a,b]$, one has,
\begin{equation}
\label{eq1haar:localTh2}
\big|\th (2^j u+k,v)\big|\le c \big (1+k\big)^{v-1/\al-1}
\end{equation}
and 
\begin{equation}
\label{eq2haar:localTh2}
\big|\Th (2^j u+k,v)\big|\le c' \big (1+k\big)^{v-1/\al-2}.
\end{equation}
\end{Prop}

The following lemma is a straightforward consequence of Lemma~\ref{lemhaar:abellf+} as well as Propositions~\ref{prophaar:estimljk2}~and~\ref{prophaar:localTh2}.

\begin{lemma}
\label{lemhaar:boundJ+A+B+}
Let $\eta>0$ be arbitrary and fixed. There exists a positive and finite random variable $C$ only depending on $\eta>0$, such that any $\o\in\Omega_{1}^{**}$, $(J,Q)\in\N^2$ 
and $v\in [a,b]$,
satisfy the following two properties:
\begin{itemize}
\item[(i)] for each $j\in\{0,\ldots,J-1\}$, one has,
\begin{align}
\label{eq1haar:boundJ+A+B+}
& \sup_{u\in [0,1]} \bigg|\sum_{k=2^{J-j}+1}^{2^{J+Q-j}} \lajnk(\o)\th (2^j u+k,v)\bigg|\\
&\le C(\o)(1+j)^{1/\al}\log^{1/\alpha+\eta}\big(2+ j\big)\Bigg\{2^{-(1-v)(J+Q-j)}(J+Q-j)^{1/\al+\eta}\nonumber\\
&\hspace{1.5cm}+2^{-(1-v)(J-j)}(J-j)^{1/\al+\eta}+\sum_{k=2^{J-j}+2}^{2^{J+Q-j}}k^{-(2-v)}\log^{1/\al+\eta}(k)\Bigg\};\nonumber
\end{align}
\item[(ii)] for each $j\in\{J,\ldots,J+Q-1\}$, one has
\begin{align}
\label{eq2haar:boundJ+A+B+}
& \sup_{u\in [0,1]} \bigg|\sum_{k=1}^{2^{J+Q-j}} \lajnk(\o)\th (2^j u+k,v)\bigg|\\
&\le C(\o)(1+j)^{1/\al}\log^{1/\alpha+\eta}\big(2+ j\big)\Bigg\{2^{-(1-v)(J+Q-j)}(J+Q-j)^{1/\al+\eta}\nonumber\\
&\hspace{6cm}+\sum_{k=2}^{2^{J+Q-j}}k^{-(2-v)}\log^{1/\al+\eta}(k)\Bigg\}.\nonumber
\end{align}
\end{itemize}
\end{lemma}

\begin{lemma}
\label{lemhaar:bound2J+A}
Let $\eta>0$ be an arbitrarily small fixed real number. There is a constant $c>0$, such that, for all $J\in\N$, for each $j\in\{0,\ldots,J-1\}$ and for any $v\in [a,b]$, one has,
\begin{equation}
\label{eq1:bound2J+A}
\sum_{k=2^{J-j}+2}^{+\infty}k^{-(2-v)}\log^{1/\al+\eta}(k)\le c\, 2^{-(1-v)(J-j)}(J-j)^{1/\al+\eta}.
\end{equation}
\end{lemma}

\begin{proof}[Proof of Lemma~\ref{lemhaar:bound2J+A}] By using the fact that $y\mapsto y^{-(2-v)}\log^{1/\al+\eta}(y)$ is a decreasing function on the interval 
$[3,+\infty)$, one has that,
\begin{equation}
\label{eq2:bound2J+A}
\sum_{k=2^{J-j}+2}^{+\infty}k^{-(2-v)}\log^{1/\al+\eta}(k)\le \int_{2^{J-j}}^{+\infty} y^{-(2-v)}\log^{1/\al+\eta}(y)\,dy.
\end{equation}
Moreover, setting in the last integral $z=2^{-(J-j)}y$, one obtains that,
\begin{eqnarray}
\label{eq3:bound2J+A}
&& \int_{2^{J-j}}^{+\infty} y^{-(2-v)}\log^{1/\al+\eta}(y)\,dy= 2^{J-j}\int_{1}^{+\infty} \big (2^{J-j}z\big)^{-(2-v)}\log^{1/\al+\eta}\big (2^{J-j}z\big)\,dz\nonumber\\
&& =2^{-(1-v)(J-j)}\int_{1}^{+\infty} z^{-(2-v)}\Big( (J-j)\log(2)+\log(z)\Big)^{1/\al+\eta}\,dz\nonumber\\
&& \le c\,2^{-(1-v)(J-j)}(J-j)^{1/\al+\eta},
\end{eqnarray}
where 
$$
c=\int_{1}^{+\infty} z^{-(2-b)}\big(1+\log(z)\big)^{1/\al+\eta}\,dz.
$$
Finally, combining (\ref{eq2:bound2J+A}) with (\ref{eq3:bound2J+A}), one gets the lemma.
\end{proof} 

\begin{lemma}
\label{lemhaar:bound3J+A}
Let $\eta>0$ be an arbitrarily small fixed real number. There exists a positive and finite random variable $C$ only depending on $b$ and $\eta$, such that any $\o\in\Omega_{1}^{**}$, $(J,Q)\in\N^2$ 
and $j\in\{0,\ldots,J-1\}$, satisfy:
\begin{align}
\label{eq1haar:bound3J+A}
& \bigg\| 2^{-j\cdot}\sum_{k=2^{J-j}+1}^{2^{J+Q-j}} \lajnk(\o)\th (2^j \cdot+k,\cdot)\bigg\|_{\ce}:=\sup_{(u,v)\in [0,1]\times [a,b]} \nonumber \\
& \hspace{6cm} \bigg| 2^{-jv}\sum_{k=2^{J-j}+1}^{2^{J+Q-j}} \lajnk(\o)\th (2^j u+k,v)\bigg|\nonumber\\
&\le C(\o)2^{-(1-b)J}J^{1/\al+\eta}\,2^{-(2a-1)j}(1+j)^{1/\al}\log^{1/\alpha+\eta}\big(2+ j\big).
\end{align}
\end{lemma}

\begin{proof}[Proof of Lemma~\ref{lemhaar:bound3J+A}] Let us set $c_1:=\sup_{n\in\N}\left\{2^{-(1-b)n}(1+n)^{1/\al+\eta}\right\}<\infty$, then one has,
\begin{equation}
\label{eq2haar:bound3J+A}
2^{-(1-v)(J+Q-j)}(J+Q-j)^{1/\al+\eta}\le c_1 2^{-(1-v)(J-j)}(J-j)^{1/\al+\eta};
\end{equation}
the inequality (\ref{eq2haar:bound3J+A}) follows from the fact that,
\begin{align*}
& 2^{-(1-v)(J+Q-j)}(J+Q-j)^{1/\al+\eta} \\
& = 2^{-(1-v)Q}\Big(1+\frac{Q}{J-j}\Big)^{1/\al+\eta}2^{-(1-v)(J-j)}(J-j)^{1/\al+\eta}\\
&\le  2^{-(1-b)Q}(1+Q)^{1/\al+\eta}\,2^{-(1-v)(J-j)}(J-j)^{1/\al+\eta}.
\end{align*}
Next putting together, (\ref{eq1haar:boundJ+A+B+}), (\ref{eq2haar:bound3J+A}) and (\ref{eq1:bound2J+A}), one gets that for any arbitrary $v\in [a,b]$,
\begin{eqnarray}
\label{eq3haar:bound3J+A}
&& \sup_{u\in [0,1]} \bigg|2^{-jv} \sum_{k=2^{J-j}+1}^{2^{J+Q-j}} \lajnk(\o)\th (2^j u+k,v)\bigg|\nonumber\\
&&\le C (\o) 2^{-(1-v)(J-j)-jv}(J-j)^{1/\al+\eta}\,(1+j)^{1/\al}\log^{1/\alpha+\eta}\big(2+ j\big)\nonumber\\
&& \le  C (\o) 2^{-(1-v)J-(2v-1)j}(J-j)^{1/\al+\eta}\,(1+j)^{1/\al}\log^{1/\alpha+\eta}\big(2+ j\big)\nonumber\\
&&\le C(\o)2^{-(1-b)J}J^{1/\al+\eta}\,2^{-(2a-1)j}\,(1+j)^{1/\al}\log^{1/\alpha+\eta}\big(2+ j\big),
\end{eqnarray}
where $C$ is a positive and finite random variable only depending on $b$ and $\eta$; thus (\ref{eq3haar:bound3J+A}) implies that (\ref{eq1haar:bound3J+A})
holds.
\end{proof}

\begin{lemma}
\label{lemhaar:bound4J+A}
Let $\eta>0$ be an arbitrarily small fixed real number. There exists a positive and finite random variable $C$ only depending on $a$, $b$ and $\eta$, such that any $\o\in\Omega_{1}^{**}$ and $(J,Q)\in\N^2$, satisfy
\begin{equation}
\label{eq1haar:bound4J+A}
\big\|A_{2,+}^{J,Q}(\cdot,\cdot,\o)\big\|_{\ce}\le C(\o)2^{-(1-b)J}J^{1/\al+\eta};
\end{equation} 
recall that $A_{2,+}^{J,Q}(\cdot,\cdot,\o)$ has been defined in (\ref{eqhaar:defA2+}).
\end{lemma}

\begin{proof}[Proof of Lemma~\ref{lemhaar:bound4J+A}] The lemma can be obtained by using (\ref{eqhaar:defA2+}), the triangle inequality, (\ref{eq1haar:bound3J+A}) 
and the fact that 
\begin{equation}
\label{eq2haar:bound4J+A}
\sum_{j=0}^{+\infty} 2^{-(2a-1)j}\,(1+j)^{1/\al}\log^{1/\alpha+\eta}\big(2+ j\big)<\infty;
\end{equation}
notice that (\ref{eq2haar:bound4J+A}) results from our assumption that $a>1/\al>1/2$.
\end{proof}

\begin{lemma}
\label{lemhaar:bound1J+B}
Let $\eta>0$ be an arbitrarily small fixed real number. There exists a positive and finite random variable $C$ only depending on $a$, $b$ and $\eta$, such that any $\o\in\Omega_{1}^{**}$ and $(J,Q)\in\N^2$, satisfy
\begin{equation}
\label{eq1haar:bound1J+B}
\big\|B_{2,+}^{J,Q}(\cdot,\cdot,\o)\big\|_{\ce}\le C(\o)\,2^{-a J}J^{1/\al}\log^{1/\alpha+\eta}\big(1+ J\big);
\end{equation} 
recall that $B_{2,+}^{J,Q}(\cdot,\cdot,\o)$ has been defined in (\ref{eqhaar:defB2+}).
\end{lemma}

\begin{proof}[Proof of Lemma~\ref{lemhaar:bound1J+B}] First notice that by using the fact that, 
$$
2^{-(1-v)(J+Q-j)}(J+Q-j)^{1/\al+\eta}+\sum_{k=2}^{2^{J+Q-j}}k^{-(2-v)}\log^{1/\al+\eta}(k)\le c_1,
$$
where the finite constant 
$$
c_1:=\Big(\sup_{n\in\N} 2^{-(1-b)n} n^{1/\al+\eta}\Big)+\sum_{k=2}^{+\infty}k^{-(2-b)}\log^{1/\al+\eta}(k);
$$
it follows from (\ref{eq2haar:boundJ+A+B+}), that for any arbitrary $v\in [a,b]$ and $j\in\{J,\ldots,J+Q-1\}$, one has 
\begin{align}
\label{eq2haar:bound1J+B}
& \sup_{u\in [0,1]} \bigg|2^{-jv} \sum_{k=1}^{2^{J+Q-j}} \lajnk(\o)\th (2^j u+k,v)\bigg| \nonumber \\
&\le  C_2(\o) 2^{-ja}(1+j)^{1/\al}\log^{1/\alpha+\eta}\big(2+ j\big),
\end{align}
where $C_2$ is a positive and finite random variable only depending on $b$ and $\eta$. Next, using (\ref{eqhaar:defB2+}), the triangle inequality,
and (\ref{eq2haar:bound1J+B}), one gets that,
\begin{equation}
\label{eq3haar:bound1J+B}
\big\|B_{2,+}^{J,Q}(\cdot,\cdot,\o)\big\|_{\ce}\le 2 C_2 (\o)\sum_{j=J}^{+\infty} 2^{-ja}(1+j)^{1/\al}\log^{1/\alpha+\eta}\big(2+ j\big);
\end{equation}
moreover, similarly to (\ref{thmain1:eq4haar}), one can show that,
\begin{equation}
\label{eq4haar:bound1J+B}
\sum_{j=J}^{+\infty} 2^{-ja}(1+j)^{1/\al}\log^{1/\alpha+\eta}\big(2+ j\big)\le c_3 2^{-Ja}(1+J)^{1/\al}\log^{1/\alpha+\eta}\big(2+ J\big),
\end{equation}
where $c_3$ is a finite constant only depending on $a$. Finally combining (\ref{eq3haar:bound1J+B}) with (\ref{eq4haar:bound1J+B}) one obtains the lemma.
\end{proof}

Now we are in position to prove Theorem~\ref{th1haar:thmain2}.

\begin{proof}[Proof of Theorem~\ref{th1haar:thmain2}] The theorem results from (\ref{eqhaar:decompX2+}), the triangle inequality,  Lemma \ref{lemhaar:bound4J+A}, Lemma \ref{lemhaar:bound1J+B}
and the inequalities $0<1-b<a$; these inequalities are consequences of our assumption that $[a,b]\subset (1/\al,1)$.
\end{proof}

\subsection{Study of the part $j<0$ of the series}
\label{subsechaar:jnegative}

\begin{thm}
\label{th2haar:thmain2}
For each integer $J\ge 2$, we denote by $X_{2,-}^J =\{X_{2,-}^J(u,v): (u,v)\in[0,1]\times [a,b]\}$ the $\sas$ field with paths in $\ce$, defined, for every 
$(u,v)\in[0,1]\times [a,b]$, as:
\begin{equation}
\label{eqhaar:partialsum-lf-}
X_{2,-}^J (u,v)=\sum_{j=1}^{J-1} 2^{jv}\sum_{k=1}^{2^{J-j}} \lanjnk\big(\th (2^{-j} u+k,v)-\th (k,v)\big).
\end{equation}
Let $\Omega^{**}_{2}$ be the event of probability 1 introduced in Proposition~\ref{prophaar:estimljk2-}. Then for all $\o\in\Omega^{**}_{2}$,  $\big(X_{2,-}^J(\cdot,\cdot,\o)\big)_{J\ge 2}$ is a Cauchy sequence in $\ce$, moreover, its limit $X_{2,-}(\cdot,\cdot,\o)$, satisfies, for 
each fixed $\eta>0$,
\begin{equation}
\label{eq1:thmain2-}
\big\|X_{2,-}(\cdot,\cdot,\o)-X_{2,-}^J(\cdot,\cdot,\o)\big\|_{\ce}=\mathcal{O} \left (2^{-J(1-b)}J^{1/\al}\log^{1/\al+\eta}\big(J\big)\right)
\end{equation}
\end{thm}

In order to show that $\big(X_{2,-}^J(\cdot,\cdot,\o)\big)_{J\ge 2}$ is a Cauchy sequence, one needs to appropriately bound the quantity $\big\|X_{2,-}^{J+Q}(\cdot,\cdot,\o)-X_{2,-}^{J}(\cdot,\cdot,\o)\big\|_{\ce}$, for all integer $J\ge 2$, and for
all $Q\ge 1$. Observe that, in view of (\ref{eqhaar:partialsum-lf-}), one has for every 
$(u,v)\in[0,1]\times [a,b]$,
\begin{equation}
\label{eqhaar:decompX2-}
X_{2,-}^{J+Q}(\cdot,\cdot,\o)-X_{2,-}^{J}(\cdot,\cdot,\o)=A_{2,-}^{J,Q}(u,v,\o)+B_{2,-}^{J,Q}(u,v,\o),
\end{equation}
where:
\begin{equation}
\label{eqhaar:defA2-}
A_{2,-}^{J,Q}(u,v,\o):=\sum_{j=1}^{J-1} 2^{jv}\sum_{k=2^{J-j}+1}^{2^{J+Q-j}} \lanjnk(\o)\big(\th (2^{-j} u+k,v)-\th (k,v)\big),
\end{equation}
and 
\begin{equation}
\label{eqhaar:defB2-}
B_{2,-}^{J,Q}(u,v,\o):=\sum_{j=J}^{J+Q-1} 2^{jv}\sum_{k=1}^{2^{J+Q-j}} \lanjnk(\o)\big(\th (2^{-j} u+k,v)-\th (k,v)\big).
\end{equation}
For all $(j,k)\in\N^2$, let $\lnjnk$ be the $\sas$ random variable of scale parameter $k^{1/\al}$ defined as:
\begin{equation}
\label{eq2haar:abellf-}
\lnjnk:=\sum_{m=1}^{k} \zeta_{-j,-m}.
\end{equation}

The proof of the following lemma is similar to that of Lemma~\ref{lemhaar:abelhf}.

\begin{lemma}
\label{lemhaar:abellf-}
Recall that the function $\Th$ has been introduced in (\ref{eq2haar:abelhf}). \\
Let $(J,Q)\in(\N\setminus\{1\})\times\N$ and $(u,v)\in [0,1]\times [a,b]$ be arbitrary and fixed.
\begin{itemize}
\item[(i)] For each $j\in\{1,\ldots,J-1\}$, one has,
\begin{align}
\label{eq3haar:abellf-}
& \sum_{k=2^{J-j}+1}^{2^{J+Q-j}} \lanjnk\big(\th (2^{-j} u+k,v)-\th (k,v)\big)\\
&=\lambda_{-j,-2^{J+Q-j}}\big(\th (2^{-j} u+2^{J+Q-j},v)-\th (2^{J+Q-j},v)\big)\nonumber\\
& -\lambda_{-j,-2^{J-j}}\big(\th (2^{-j} u+2^{J-j}+1,v)-\th (2^{J-j}+1,v)\big)\nonumber\\
& -\sum_{k=2^{J-j}+2}^{2^{J+Q-j}}\lambda_{-j,-(k-1)}\big(\Th(2^{-j} u+k,v)-\Th(k,v)\big).\nonumber
\end{align}
\item[(ii)] For each $j\in\{J,\ldots,J+Q-1\}$, one has
\begin{align}
\label{eq4haar:abellf-}
& \sum_{k=1}^{2^{J+Q-j}} \lanjnk\big(\th (2^{-j} u+k,v)-\th (k,v)\big)\\
& =\lambda_{-j,-2^{J+Q-j}}\big(\th (2^{-j} u+2^{J+Q-j},v)-\th (2^{J+Q-j},v)\big)\nonumber\\
&-\sum_{k=2}^{2^{J+Q-j}}\lambda_{-j,-(k-1)}\big (\Th(2^{-j} u+k,v)-\Th (k,v)\big).\nonumber
\end{align}
\end{itemize}
\end{lemma}

The proof of the following proposition is similar to that of Proposition~\ref{prophaar:estimljk}.

\begin{Prop}
\label{prophaar:estimljk2-}
There exists an event of probability 1, denoted by $\Omega_2^{**}$, such that for every fixed real number $\eta >0$, one has, for all $\omega \in \Omega_2^{**}$ and for each $(j,k)\in\Z_{+}\times\N$,
\begin{equation}
\label{eq1:haar:estimljk2-}
\big| \lnjnk (\o)\big| \leq C(\omega) \big(1+j \big)^{1/\alpha}\log^{1/\alpha+\eta}\big(3+ j\big)\,  k^{1/\alpha} \log^{1/\alpha+\eta}\big(2+ k\big), 
\end{equation}
where $C$ is a positive and finite random variable only depending on $\eta$.
\end{Prop}

We denote by $\partial_{x}\th$ and $\partial_{x}\Th$, the partial derivatives of order $1$ with respect to $x$ of the functions $\th$ and $\Th$; observe that in view of (\ref{theta:00uv:alpha}), (\ref{eqhaar:adpat+}) and (\ref{eq2haar:abelhf}), one has for all $(x,v)\in\R\times (1/\al,1)$,
\begin{equation}
\label{eqhaar:partialxth}
(\partial_{x}\th)(x,v)=(x-1)_{+}^{v-1/\al}-2(x-1/2)_{+}^{v-1/\al}+ (x)_{+}^{v-1/\al}
\end{equation}
and 
\begin{equation}
\label{eqhaar:partialxTh}
(\partial_{x}\Th)(x,v)=(\partial_{x}\th)(x,v)-(\partial_{x}\th)(x-1,v)=\sum_{l=0}^{4} d_l (x-l/2)_{+}^{v-1/\al},
\end{equation}
where $d_0,\ldots,d_4$ have been defined in (\ref{eqhaar:defdl}). The proof of the following proposition relies on (\ref{eqhaar:partialxth}) and (\ref{eqhaar:partialxTh}); we will not give it since it is very similar to that of Proposition~\ref{prophaar:localTh}.

\begin{Prop}
\label{prophaar:localpartialxTh}
\begin{itemize}
\item[(i)] For each $(x,v)\in (-\infty,0]\times (1/\al,1)$, one has, 
$$(\partial_{x}\th)(x,v)=(\partial_{x}\Th)(x,v)=0.$$
\item[(ii)] There exists a constant $c>0$, such that for all $(x,v)\in (0,+\infty)\times [a,b]$, one has,
$$
|(\partial_{x}\th)(x,v)|\le c(1+x)^{v-1/\al-2}.
$$
\item[(iii)] There is a constant $c'>0$, such that for all $(x,v)\in (0,+\infty)\times [a,b]$, one has,
$$
|(\partial_{x}\Th)(x,v)|\le c'(1+x)^{v-1/\al-3}.
$$
\end{itemize} 
\end{Prop}

The following proposition easily results from the Mean Value Theorem and from Proposition~\ref{prophaar:localpartialxTh}.
\begin{Prop}
\label{prophaar:localTh2-}
There exist two constants $c>0$ and $c'>0$, such that for all $(j,k)\in\N^2$ and $(u,v)\in [0,1]\times [a,b]$, one has,
\begin{equation}
\label{eq1haar:localTh2-}
\big|\th (2^{-j} u+k,v)-\th (k,v)\big|\le c 2^{-j}\big (1+k\big)^{v-1/\al-2}
\end{equation}
and 
\begin{equation}
\label{eq2haar:localTh2-}
\big|\Th (2^{-j} u+k,v)-\Th (k,v)\big|\le c'  2^{-j}\big (1+k\big)^{v-1/\al-3}.
\end{equation}
\end{Prop}

The following lemma is a straightforward consequence of Lemma~\ref{lemhaar:abellf-} as well as Propositions~\ref{prophaar:estimljk2-}~and~\ref{prophaar:localTh2-}.

\begin{lemma}
\label{lemhaar:boundJ-A-B-}
Let $\eta>0$ be arbitrary and fixed. There exists a positive and finite random variable $C$ only depending on $\eta>0$, such that any $\o\in\Omega_{2}^{**}$, $(J,Q)\in(\N\setminus\{1\})\times\N$ 
and $v\in [a,b]$,
satisfy the following two properties:
\begin{itemize}
\item[(i)] for each $j\in\{1,\ldots,J-1\}$, one has,
\begin{align}
\label{eq1haar:boundJ-A-B-}
& \sup_{u\in [0,1]} \bigg|\sum_{k=2^{J-j}+1}^{2^{J+Q-j}} \lanjnk(\o)\big(\th (2^{-j} u+k,v)-\th (k,v)\big)\bigg| \nonumber \\
&\le C(\o)2^{-j}(1+j)^{1/\al}\log^{1/\alpha+\eta}\big(2+ j\big)\Bigg\{2^{-(2-v)(J+Q-j)}(J+Q-j)^{1/\al+\eta}\nonumber\\
&\hspace{2cm}+2^{-(2-v)(J-j)}(J-j)^{1/\al+\eta}+\sum_{k=2^{J-j}+2}^{2^{J+Q-j}}k^{-(3-v)}\log^{1/\al+\eta}(k)\Bigg\};
\end{align}
\item[(ii)] for each $j\in\{J,\ldots,J+Q-1\}$, one has
\begin{align}
\label{eq2haar:boundJ-A-B-}
& \sup_{u\in [0,1]} \bigg|\sum_{k=1}^{2^{J+Q-j}} \lanjnk(\o)\big(\th (2^{-j} u+k,v)-\th (k,v)\big)\bigg| \nonumber \\
&\le C(\o)2^{-j}(1+j)^{1/\al}\log^{1/\alpha+\eta}\big(2+ j\big)\Bigg\{2^{-(2-v)(J+Q-j)}(J+Q-j)^{1/\al+\eta}\nonumber\\
&\hspace{5cm}+\sum_{k=2}^{2^{J+Q-j}}k^{-(3-v)}\log^{1/\al+\eta}(k)\Bigg\}.
\end{align}
\end{itemize}
\end{lemma}

The proof of the following lemma mainly relies on (\ref{eq1haar:boundJ-A-B-}), it can be done similarly to that of Lemma~\ref{lemhaar:bound3J+A}.
\begin{lemma}
\label{lemhaar:bound3J-A}
Let $\eta>0$ be an arbitrarily small fixed real number. There exists a positive and finite random variable $C$ only depending on $b$ and $\eta$, such that any $\o\in\Omega_{2}^{**}$, $(J,Q)\in(\N\setminus\{1\})\times\N$ 
and $j\in\{1,\ldots,J-1\}$, satisfy:
\begin{eqnarray}
\label{eq1haar:bound3J-A}
&& \bigg\| 2^{j\cdot}\sum_{k=2^{J-j}+1}^{2^{J+Q-j}} \lanjnk(\o)\big(\th (2^{-j} \cdot+k,\cdot)-\th (k,\cdot)\big)\bigg\|_{\ce}\nonumber\\
&&:=\sup_{(u,v)\in [0,1]\times [a,b]}\bigg| 2^{jv}\sum_{k=2^{J-j}+1}^{2^{J+Q-j}} \lanjnk(\o)\big(\th (2^{-j} \cdot+k,\cdot)-\th (k,\cdot)\big)\bigg|\nonumber\\
&&\le C(\o)2^{-(1-b)J}2^{-(J-j)}(J-j)^{1/\al+\eta}\,(1+j)^{1/\al}\log^{1/\alpha+\eta}\big(2+ j\big).
\end{eqnarray}
\end{lemma}

\begin{lemma}
\label{lemhaar:bound4J-A}
Let $\eta>0$ be an arbitrarily small fixed real number. There exists a positive and finite random variable $C$ only depending on $b$ and $\eta$, such that any $\o\in\Omega_{2}^{**}$ and $(J,Q)\in(\N\setminus\{1\})\times\N$, satisfy
\begin{equation}
\label{eq1haar:bound4J-A}
\big\|A_{2,-}^{J,Q}(\cdot,\cdot,\o)\big\|_{\ce}\le C(\o)2^{-(1-b)J}J^{1/\al}\log^{1/\alpha+\eta}\big(1+J\big);
\end{equation} 
recall that $A_{2,-}^{J,Q}(\cdot,\cdot,\o)$ has been defined in (\ref{eqhaar:defA2-}).
\end{lemma}

\begin{proof}[Proof of Lemma~\ref{lemhaar:bound4J-A}] Using (\ref{eqhaar:defA2+}), the triangle inequality, and (\ref{eq1haar:bound3J+A}), one 
obtains that, 
\begin{eqnarray}
\label{eq2haar:bound4J-A}
&& \big\|A_{2,-}^{J,Q}(\cdot,\cdot,\o)\big\|_{\ce}\\
&& \le  C(\o)2^{-(1-b)J}\sum_{j=1}^{J-1} 2^{-(J-j)}(J-j)^{1/\al+\eta}\,(1+j)^{1/\al}\log^{1/\alpha+\eta}\big(2+ j\big)\nonumber\\
&& \le  C(\o)2^{-(1-b)J}J^{1/\al}\log^{1/\alpha+\eta}\big(1+ J\big)\sum_{j=1}^{J-1} 2^{-(J-j)}(J-j)^{1/\al+\eta};\nonumber
\end{eqnarray}
moreover, observe that,
\begin{equation}
\label{eq3haar:bound4J-A}
\sum_{j=1}^{J-1} 2^{-(J-j)}(J-j)^{1/\al+\eta}\le c_1,
\end{equation}
where $c_1$ is the finite constant non depending on $J$, defined as 
$$
c_1=\sum_{n=1}^{+\infty} 2^{-n} n^{1/\al+\eta}.
$$
Finally, combining (\ref{eq2haar:bound4J-A}) with (\ref{eq3haar:bound4J-A}) one gets the lemmma.
\end{proof}

\begin{lemma}
\label{lemhaar:bound1J-B}
Let $\eta>0$ be an arbitrarily small fixed real number. There exists a positive and finite random variable $C$ only depending on $b$ and $\eta$, such that any $\o\in\Omega_{2}^{**}$ and $(J,Q)\in(\N\setminus\{1\})\times\N$, satisfy
\begin{equation}
\label{eq1haar:bound1J-B}
\big\|B_{2,-}^{J,Q}(\cdot,\cdot,\o)\big\|_{\ce}\le C(\o)\,2^{-(1-b) J}J^{1/\al}\log^{1/\alpha+\eta}\big(1+ J\big);
\end{equation} 
recall that $B_{2,-}^{J,Q}(\cdot,\cdot,\o)$ has been defined in (\ref{eqhaar:defB2-}).
\end{lemma}

\begin{proof}[Proof of Lemma~\ref{lemhaar:bound1J-B}] First notice that by using the fact that, 
$$
2^{-(2-v)(J+Q-j)}(J+Q-j)^{1/\al+\eta}+\sum_{k=2}^{2^{J+Q-j}}k^{-(3-v)}\log^{1/\al+\eta}(k)\le c_1,
$$
where the finite constant 
$$
c_1:=\Big(\sup_{n\in\N} 2^{-(2-b)n} n^{1/\al+\eta}\Big)+\sum_{k=2}^{+\infty}k^{-(3-b)}\log^{1/\al+\eta}(k);
$$
it follows from (\ref{eq2haar:boundJ-A-B-}), that for any arbitrary $v\in [a,b]$ and $j\in\{J,\ldots,J+Q-1\}$, one has 
\begin{align}
\label{eq2haar:bound1J-B}
& \sup_{u\in [0,1]} \bigg|2^{jv} \sum_{k=1}^{2^{J+Q-j}} \lanjnk(\o)\big(\th (2^{-j} u+k,v)-\th (k,v)\big)\bigg| \nonumber \\
& \le  C_2(\o) 2^{-j(1-b)}(1+j)^{1/\al}\log^{1/\alpha+\eta}\big(2+ j\big),
\end{align}
where, $C_2$ is a positive and finite random variable only depending on $b$ and $\eta$. Next, using (\ref{eqhaar:defB2-}), the triangle inequality,
and (\ref{eq2haar:bound1J-B}), one gets that,
\begin{equation}
\label{eq3haar:bound1J-B}
\big\|B_{2,-}^{J,Q}(\cdot,\cdot,\o)\big\|_{\ce}\le C_2 (\o)\sum_{j=J}^{+\infty} 2^{-j(1-b)}(1+j)^{1/\al}\log^{1/\alpha+\eta}\big(2+ j\big);
\end{equation}
Moreover, similarly to (\ref{thmain1:eq4haar}), one can show that,
\begin{equation}
\label{eq4haar:bound1J-B}
\sum_{j=J}^{+\infty} 2^{-j(1-b)}(1+j)^{1/\al}\log^{1/\alpha+\eta}\big(2+ j\big)\le c_3 2^{-J(1-b)}(1+J)^{1/\al}\log^{1/\alpha+\eta}\big(2+ J\big),
\end{equation}
where $c_3$ is a finite constant only depending on $b$. Finally combining (\ref{eq3haar:bound1J-B}) with (\ref{eq4haar:bound1J-B}) one obtains the lemma.
\end{proof}

Now we are in position to prove Theorems~\ref{th2haar:thmain2}~and~\ref{thhaar:thmain2}.

\begin{proof}[Proof of Theorem~\ref{th2haar:thmain2}] The theorem results from (\ref{eqhaar:decompX2-}) as well as  Lemmas \ref{lemhaar:bound4J-A} and \ref{lemhaar:bound1J-B}.
\end{proof}

\begin{proof}[Proof of Theorem~\ref{thhaar:thmain2}] Putting together (\ref{eqhaar:partialsum-lf}), (\ref{eqhaar:partialsum-lf+}), 
(\ref{eqhaar:partialsum-lf-}), Theorem~\ref{th1haar:thmain2} and Theorem~\ref{th2haar:thmain2}, one gets the theorem.

\end{proof}

\section{Simulations}\label{sechaar:simulation}

Let us stress that Theorem \ref{thhaar:thmain1} and Theorem \ref{thhaar:thmain2} provide an efficient method for simulating paths of the high frequency part and the low frequency part of LMSM, namely of the $\sas$ processes 
$\{Y_1(t):t\in [0,1]\}:=\{X_1 (t,H(t)):t\in [0,1]\}$ and $\{Y_2(t):t\in [0,1]\}:=\{X_2 (t,H(t)):t\in [0,1]\}$. The following four simulations have been performed by using (\ref{eqhaar:partialsum-hf}) in which $J=12$ and (\ref{eqhaar:partialsum-lf}) in which $J=6$.

\begin{figure}
\begin{center}
\begin{tabular}{ccc}
\includegraphics[scale=0.22]{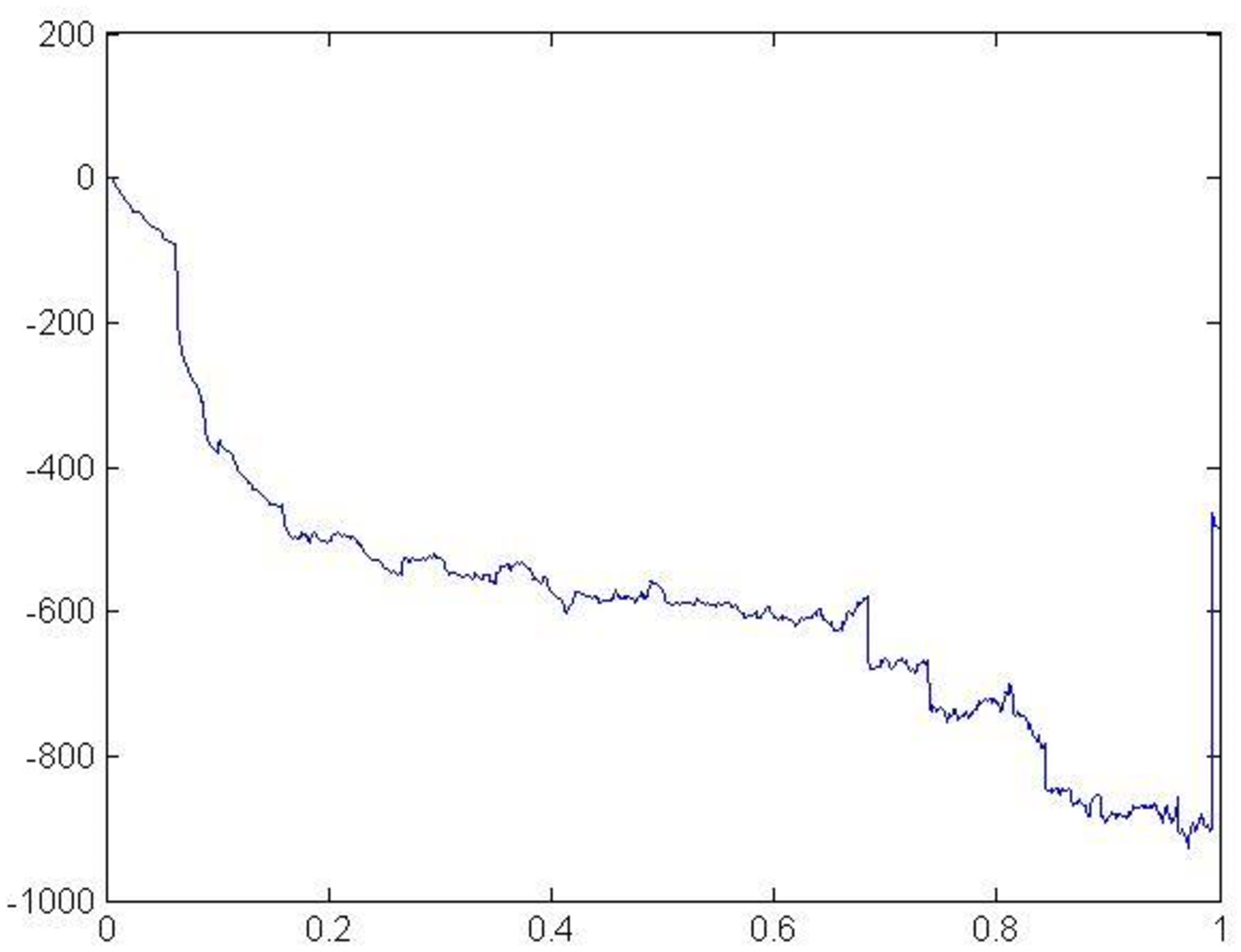} & \includegraphics[scale=0.22]{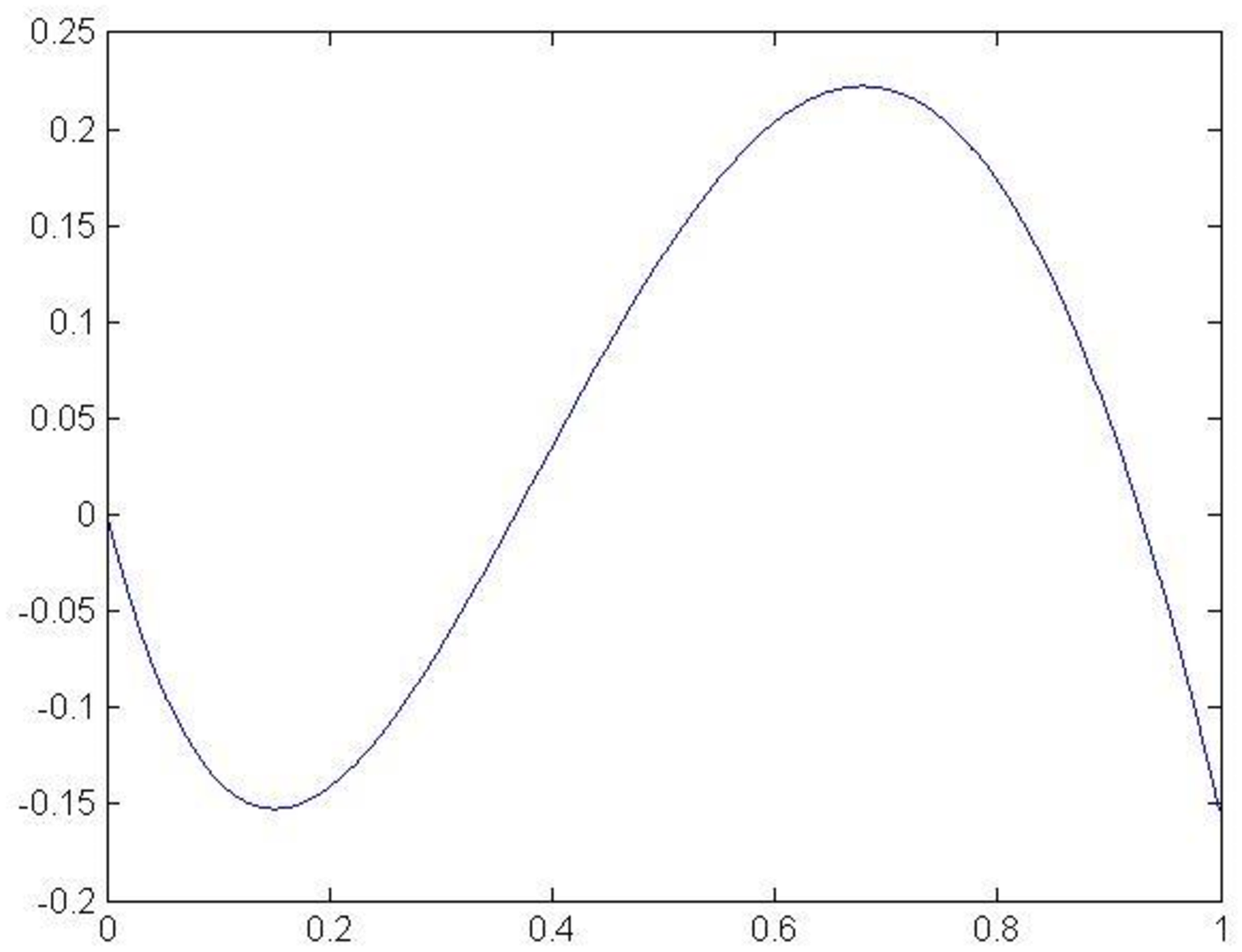} &
\includegraphics[scale=0.22]{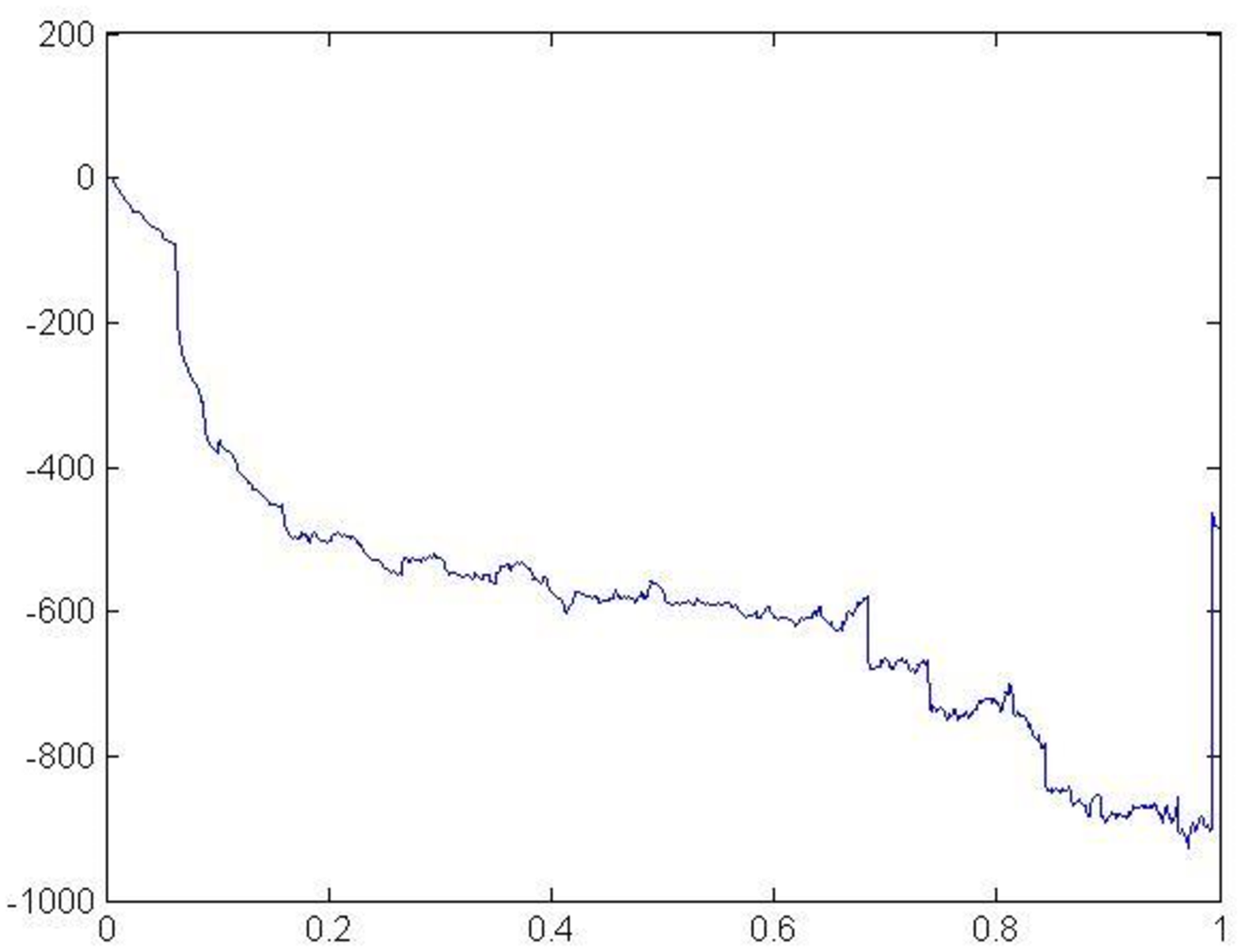} \\
\includegraphics[scale=0.22]{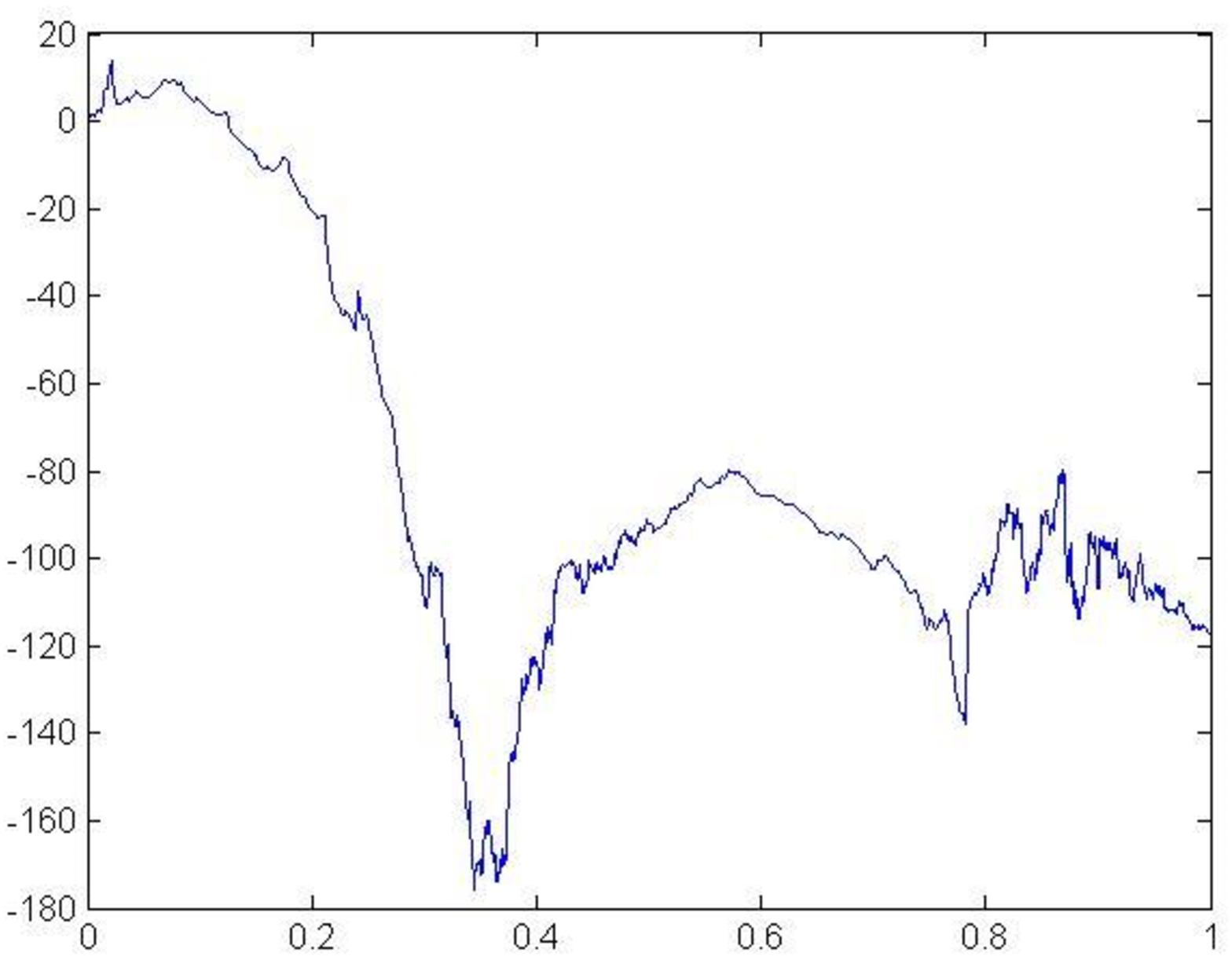} &\includegraphics[scale=0.22]{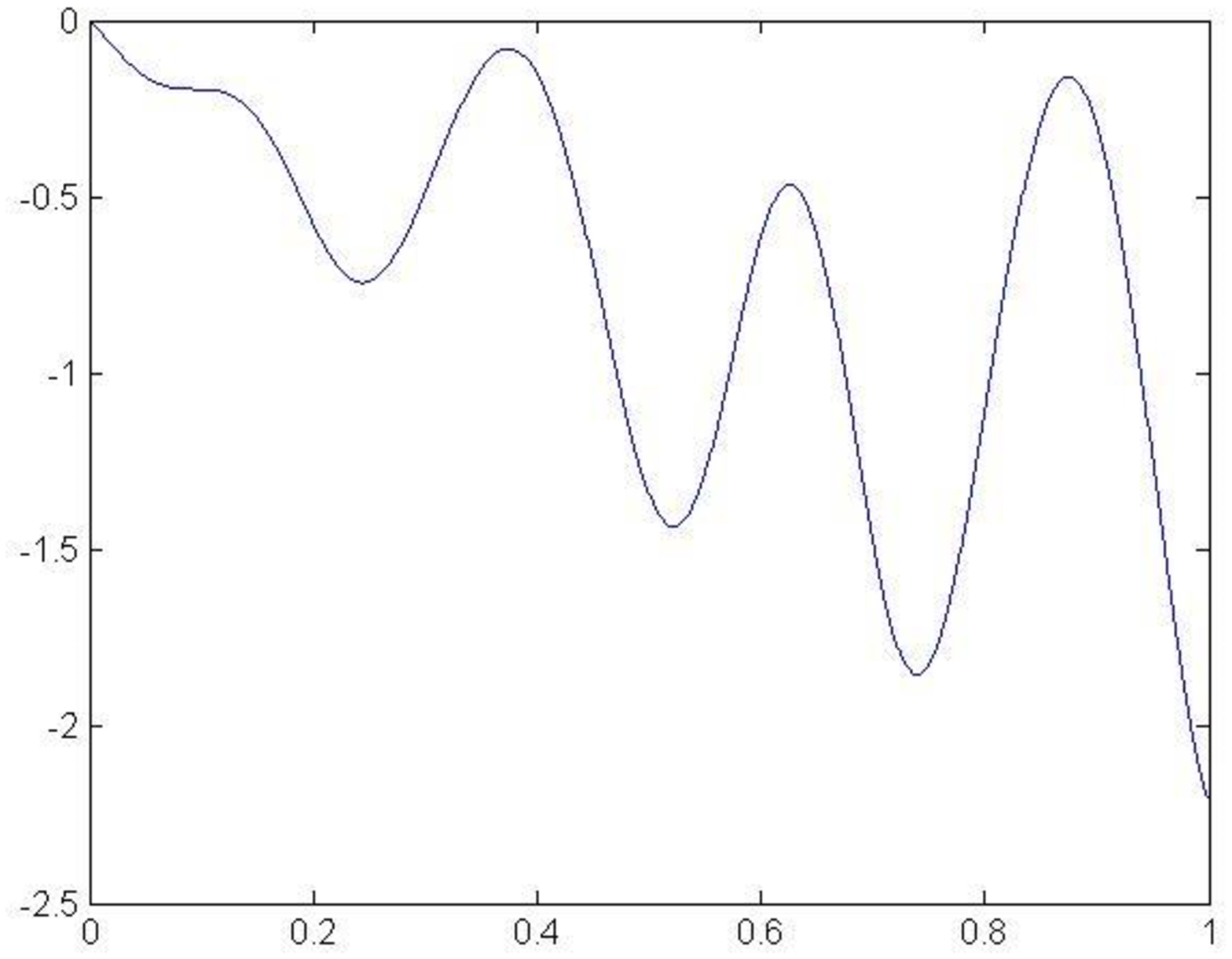} &
\includegraphics[scale=0.22]{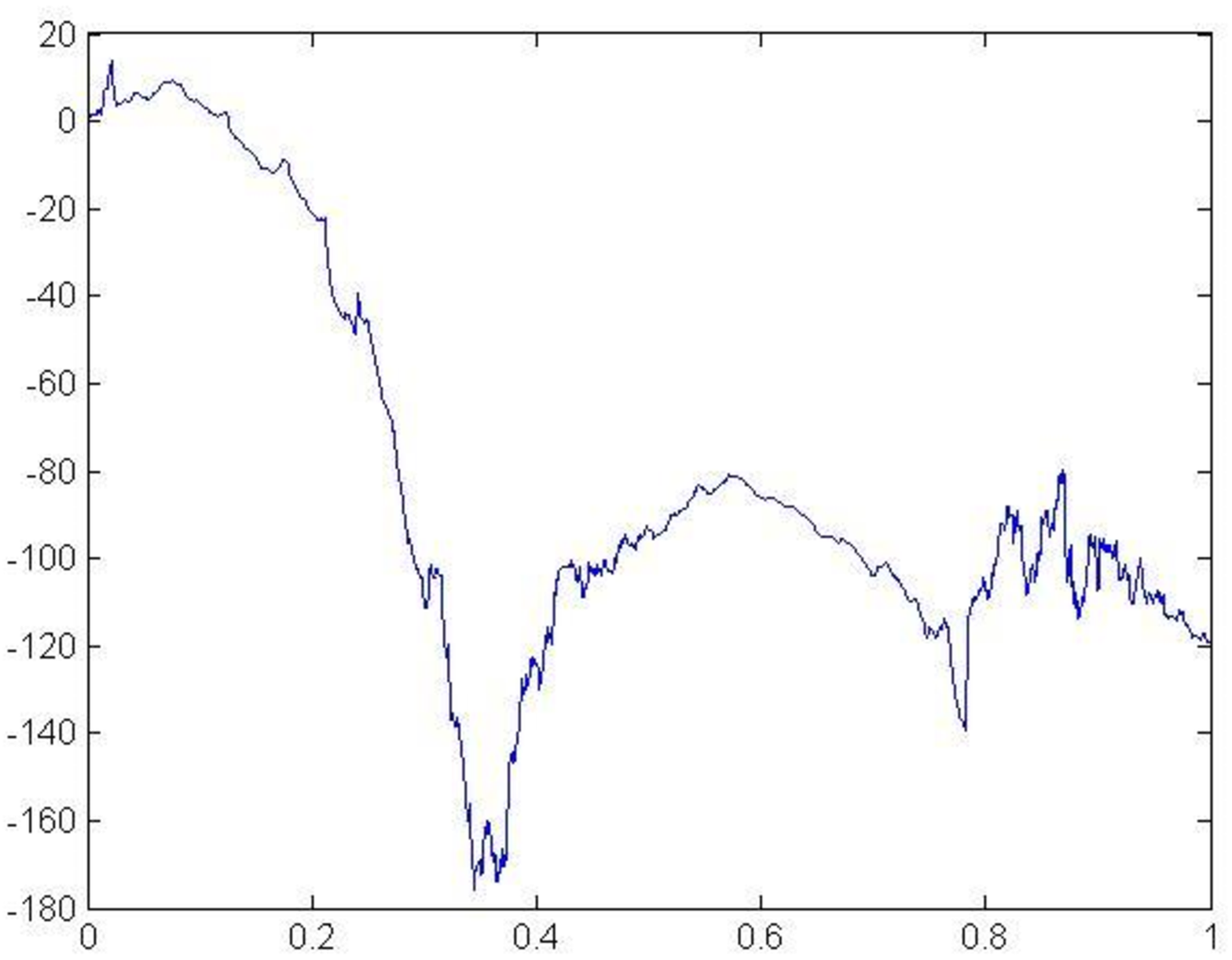} \\
\includegraphics[scale=0.22]{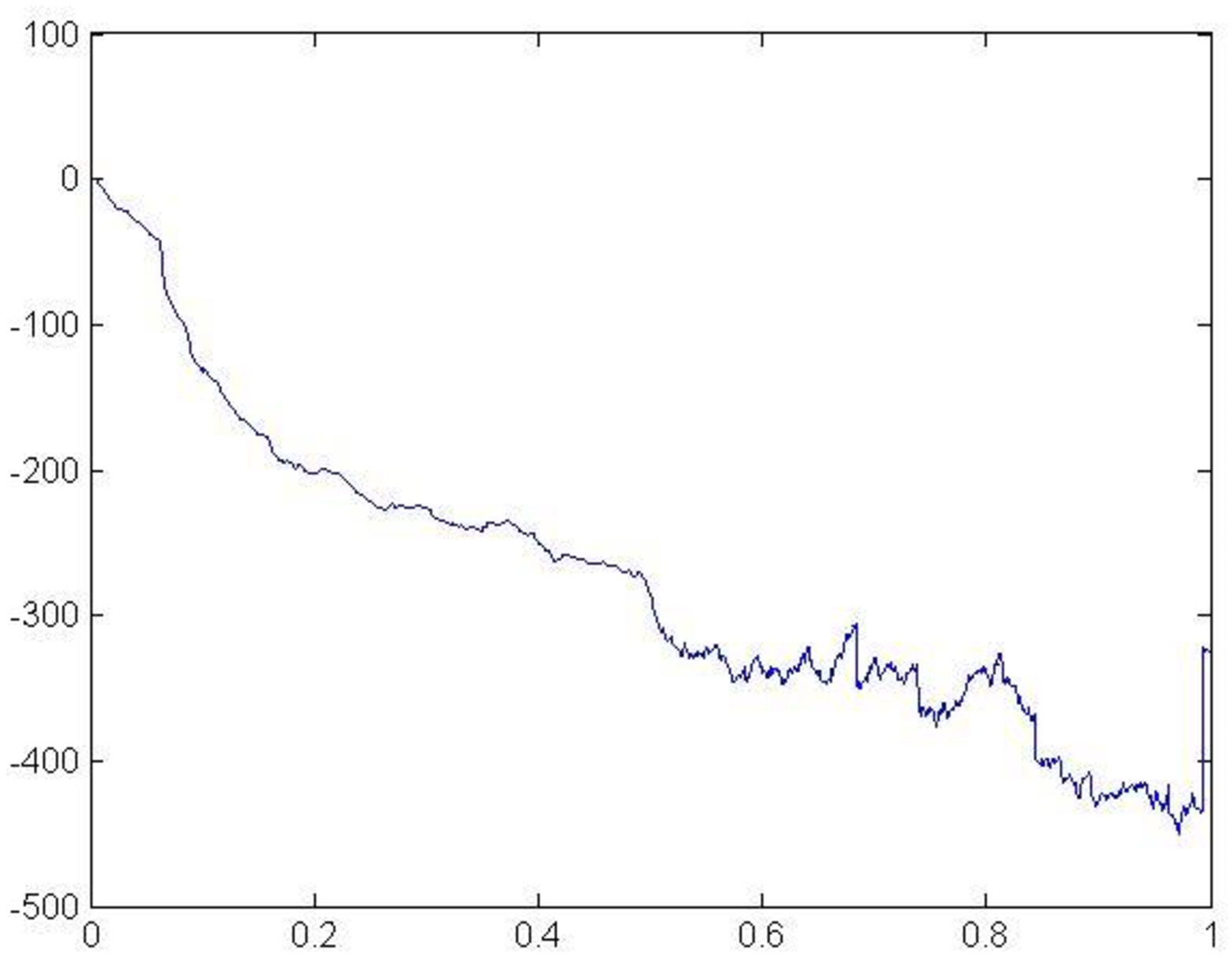} & \includegraphics[scale=0.22]{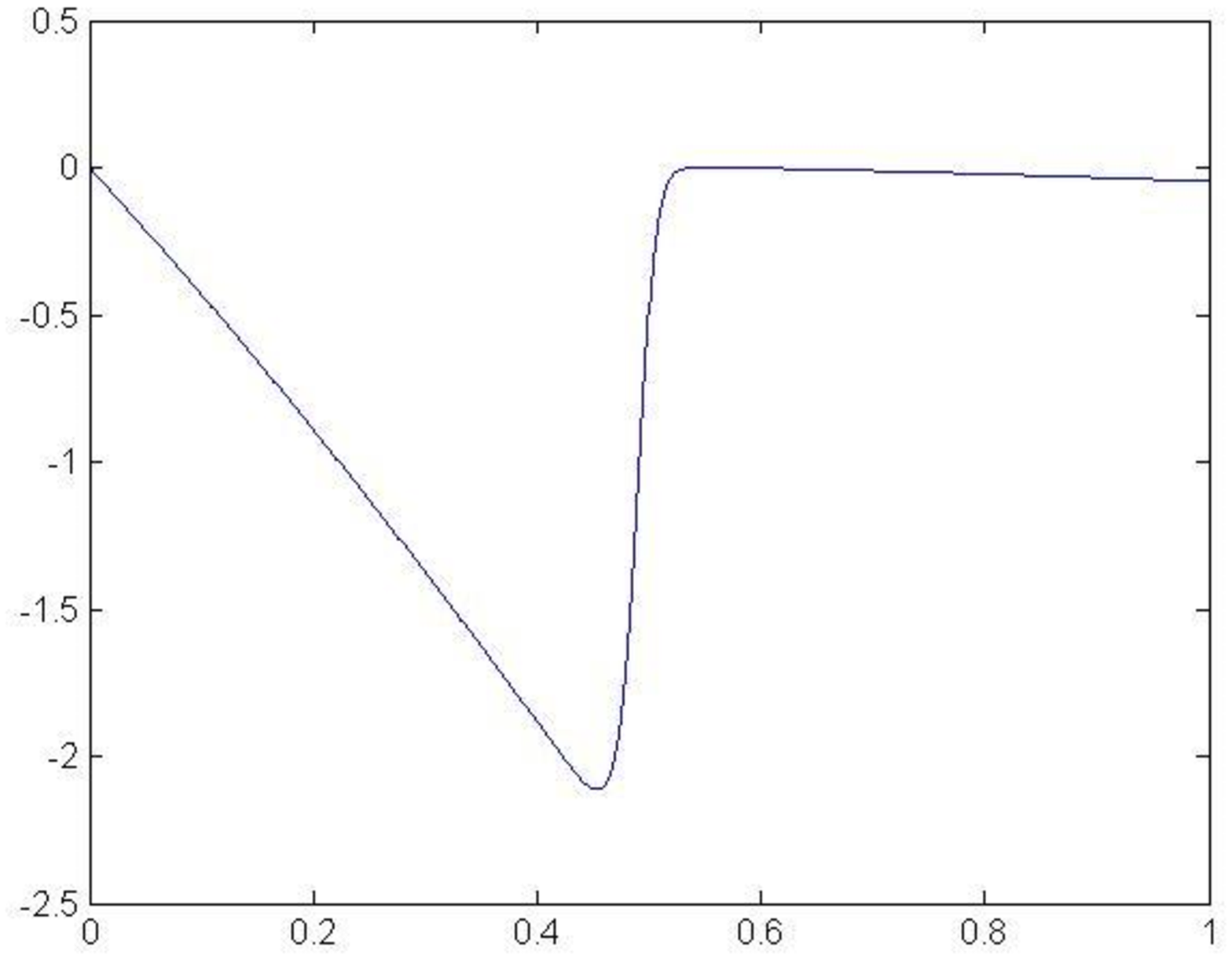} &
\includegraphics[scale=0.22]{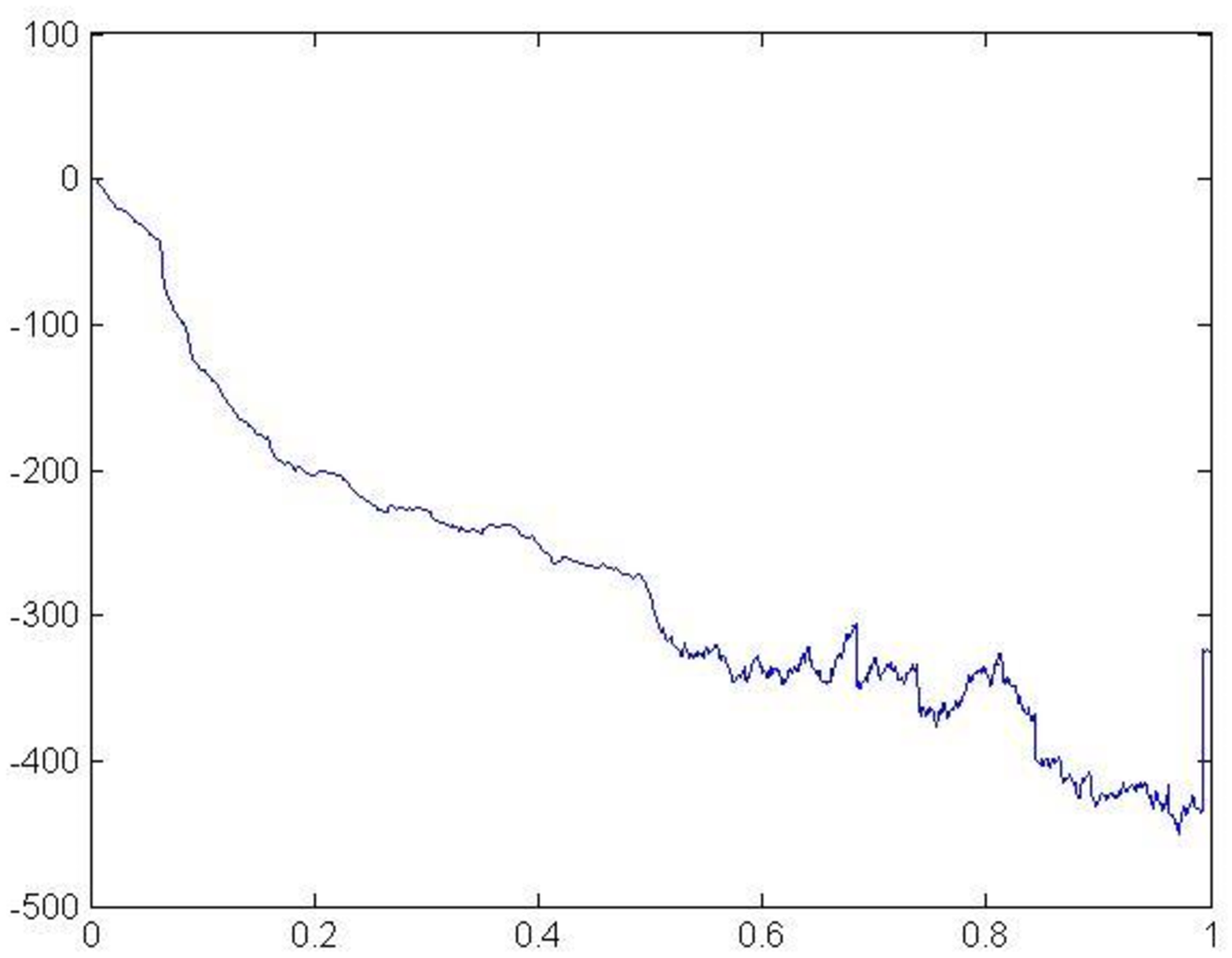}
\end{tabular}
\end{center}
\caption{In the first line: simulation of a path of the processes $Y_1$, $Y_2$ and $Y$ when $\al=1.4$ and $H(t)=0.9-0.2t$ for all $t\in [0,1]$. In the second line: simulation of a path of the processes $Y_1$, $Y_2$ and $Y$ when $\al=1.7$ and $H(t)=0.2\sin(4\pi t)+0.8$ for all $t\in [0,1]$. In the third line: simulation of a path of the processes $Y_1$, $Y_2$ and $Y$ when $\al=1.6$ and $H(t)=0.65+0.25/(1+\exp(100(t-0.5)))$ for all $t\in [0,1]$.
}
\end{figure}

\noindent{\bf Acknowlegment}
\\

I thank Professor A. Ayache for several very helpful discussions on the subject of the paper.





\section*{References}
\bibliographystyle{elsart-num-sort}
\bibliography{haarbiblio}







\end{document}